\documentclass[11pt]{amsart}

\usepackage{tikz,varwidth, hyperref, amssymb, mathscinet}
\usepackage[left=0.97in,top=0.97in,right=0.97in,bottom=0.97in]{geometry}

\usetikzlibrary{calc, snakes}
\usepackage{geometry, framed}
\usepackage{stmaryrd}
\usepackage[shortlabels]{enumitem}
\usepackage[all,cmtip]{xy}
\usepackage{tikz-cd}
\usepackage{diagbox}

\newtheorem{thm}{Theorem}[section]

\theoremstyle{definition} % text is not in italics

\newtheorem{innercustomthm}{Theorem}
\newenvironment{customthm}[1]
  {\renewcommand\theinnercustomthm{#1}\innercustomthm}
  {\endinnercustomthm}

\newtheorem{innercustomprop}{Proposition}
\newenvironment{customprop}[1]
  {\renewcommand\theinnercustomprop{#1}\innercustomprop}
  {\endinnercustomprop}
  
\usetikzlibrary{arrows}
\usetikzlibrary{decorations.markings,decorations.pathmorphing,arrows}
\pgfarrowsdeclare{mytip}{mytip}{%
	\setlength{\arrowsize}{1pt}
	\addtolength{\arrowsize}{.5\pgflinewidth}
	\pgfarrowsrightextend{0}
	\pgfarrowsleftextend{-5\arrowsize}
}{%
	\setlength{\arrowsize}{1pt}
	\addtolength{\arrowsize}{.5\pgflinewidth}
	\pgfpathmoveto{\pgfpoint{-5\arrowsize}{1.5\arrowsize}}
	\pgfpathlineto{\pgfpointorigin}
	\pgfpathlineto{\pgfpoint{-5\arrowsize}{-1.5\arrowsize}}
	\pgfusepathqstroke
}
\tikzset{myarrow/.style={ decoration={bent,aspect=0.3, markings,mark=at
			position 0.5 with {\arrow[scale=1.2]{latex'}}}, postaction=decorate}}
\tikzset{myarrowshort/.style={ decoration={bent,aspect=0.3, markings,mark=at
			position 0.3 with {\arrow[scale=1.2]{latex'}}}, postaction=decorate}}
\tikzset{myarrowshorter/.style={ decoration={bent,aspect=0.3, markings,mark=at
			position 0.2 with {\arrow[scale=1.2]{latex'}}}, postaction=decorate}}
\tikzset{->-/.style={decoration={markings, mark=at position #1 with
			{\arrow{>}}},postaction={decorate}}}

\tikzset{my_dot/.style={fill, circle, inner sep=0pt,minimum size=1.5pt}}
\tikzset{my_node/.style={fill, circle, inner sep=0pt,minimum size=3pt}}
\tikzset{inv/.style={fill, circle, inner sep=0pt,minimum size=0pt}}

\newcommand{\CC}{\mathbb{C}}

\newcommand{\FF}{\mathbb{F}}

\newcommand{\PP}{\mathbb{P}}
\newcommand{\QQ}{\mathbb{Q}}
\newcommand{\RR}{\mathbb{R}}

\newcommand{\ZZ}{\mathbb{Z}}

\newcommand{\cB}{\mathcal{B}}
\newcommand{\cC}{\mathcal{C}}

\newcommand{\cG}{\mathcal{G}}
\newcommand{\cH}{\mathcal{H}}

\newcommand{\cL}{\mathcal{L}}
\newcommand{\cM}{\mathcal{M}}

\newcommand{\cP}{\mathcal{P}}
\newcommand{\cQ}{\mathcal{Q}}
\newcommand{\cR}{\mathcal{R}}
\newcommand{\cS}{\mathcal{S}}

\newcommand{\comment}[1]{}
\newcommand{\on}{\operatorname}

\newcommand{\Spec}{\operatorname{Spec}}

\newcommand{\ov}{\overline}

\newcommand{\Aut}{\on{Aut}}

\newcommand{\Conf}{\on{Conf}}

\newcommand{\xra}{\xrightarrow}

\newcommand{\wt}{\widetilde}

\newcommand{\ch}{\on{ch}}

\newcommand{\Hom}{\on{Hom}}

\newcommand{\Adm}{\mathrm{Adm}}

\renewcommand{\P}{\mathbf{P}}
\newcommand{\C}{\mathbf{C}}

\newcommand{\col}{\colon}

\newcommand{\hide}[1]{}
\newcommand{\stub}[1]{}

\newcommand{\down}[2]{\xymatrix@R=6mm@C=2mm{
#1\ar[d]\\ #2
}}

\newcommand{\downlabel}[3]{\xymatrix@R=6mm@C=2mm{
{#1}\ar[d]^<<<{#3} \\ #2
}}

\newcommand{\squarediagram}[4]{\xymatrix@R=8mm@C=8mm{
#1\ar[d]\ar[r] & #2\ar[d] \\ #3\ar[r] &#4
}}
\newcommand{\squarediagrammapsto}[4]{\xymatrix@R=8mm@C=8mm{
#1\ar@{|->}[d]\ar@{|->}[r] & #2\ar@{|->}[d] \\ #3\ar@{|->}[r] &#4
}}
\newcommand{\squarediagramlabel}[8]{\xymatrix@R=8mm@C=8mm{
#1\ar[d]_{#6}\ar[r]^{#5} & #2\ar[d]^{#7} \\ #3\ar[r]^{#8} &#4
}}

\newcommand{\isocelesdown}[3]{\xymatrix@R=6mm@C=0mm{
& {#1}\ar[dl] \ar[dr] & \\
{#2} \ar[rr] && {#3}
}}

\newcommand{\isocelesdownlabel}[6]{\xymatrix@R=6mm@C=0mm{
& {#1}\ar[dl]_<<<<{#4} \ar[dr]^<<<<{#5} & \\
{#2} \ar[rr]_{#6} && {#3}
}}

\newcommand{\isocelesup}[3]{\xymatrix@R=6mm@C=0mm{
 #1\ar[rr]\ar[dr]  && #2\ar[dl] \\
 & #3 &
}}

\newcommand{\isocelesuplabel}[6]{\xymatrix@R=6mm@C=0mm{
 #1\ar[rr]^{{#4}} \ar[dr]_<<<{#5} && #2\ar[dl]^<<<{#6} \\
 & #3 &
}}

\newtheorem{Definition}[thm]{Definition}
\newenvironment{definition}
  {\begin{Definition}\rm}{\end{Definition}}

\newtheorem{Example}[thm]{Example}
\newenvironment{example}
  {\begin{Example}\rm}{\end{Example}}
  
\newtheorem{Exercise}[thm]{Exercise}
\newenvironment{exercise}
  {\begin{Exercise}\rm}{\end{Exercise}}
    
\newtheorem{Exploration}[thm]{Exploration}
\newenvironment{exploration}
  {\begin{Exploration}\rm}{\end{Exploration}}

\newtheorem{Fact}[thm]{Fact}
\newenvironment{fact}
  {\begin{Fact}\rm}{\end{Fact}}

\newtheorem{Theorem}[thm]{Theorem}
\newenvironment{theorem}
  {\begin{Theorem}\rm}{\end{Theorem}}
  
\newtheorem{Lemma}[thm]{Lemma}
\newenvironment{lemma}
  {\begin{Lemma}\rm}{\end{Lemma}}

\newtheorem{Remark}[thm]{Remark}
\newenvironment{remark}
  {\begin{Remark}\rm}{\end{Remark}}

\newtheorem{Proposition}[thm]{Proposition}
\newenvironment{proposition}
  {\begin{Proposition}\rm}{\end{Proposition}}

\newtheorem{Corollary}[thm]{Corollary}
\newenvironment{corollary}
  {\begin{Corollary}\rm}{\end{Corollary}}

\newtheorem{Question}[thm]{Question}
\newenvironment{question}
  {\begin{Question}\rm}{\end{Question}}

\newtheorem{Conjecture}[thm]{Conjecture}
\newenvironment{conjecture}
  {\begin{Conjecture}\rm}{\end{Conjecture}}
  
  \newtheorem{Problem}[thm]{Problem}

\theoremstyle{remark}

\newcommand \defnow[1]{\begin{definition}{#1}\end{definition}}

\newcommand \thmnow[1]{\begin{theorem}{#1}\end{theorem}}
\newcommand \proofnow[1]{\begin{proof}{#1}\end{proof}}
\newcommand \remnow[1]{\begin{remark}{#1}\end{remark}}
\newcommand \propnow[1]{\begin{proposition}{#1}\end{proposition}}

\newcommand \enumnow[1]{\begin{enumerate}{#1}\end{enumerate}}

\title{On the weight zero compactly supported cohomology of $\cH_{g,n}$}
\author[M. Brandt]{Madeline Brandt}\address{Department of Mathematics, Brown University, Box
1917, Providence, RI 02912}\email{madeline\_brandt@brown.edu}
\author[M. Chan]{Melody Chan}\address{Department of Mathematics, Brown University, Box
1917, Providence, RI 02912}\email{melody\_chan@brown.edu}
\author[S. Kannan]{Siddarth Kannan}\address{Department of Mathematics, Brown University, Box
1917, Providence, RI 02912}\email{siddarth\_kannan@brown.edu}
\date{\today}

\begin{document}
\maketitle

%\tableofcontents

\begin{abstract}
    For $g\ge 2$ and $n\ge 0$, let $\cH_{g,n}\subset \cM_{g,n}$ denote the complex moduli stack of $n$-marked smooth hyperelliptic curves of genus $g$. A normal crossings compactification of this space is provided by the theory of pointed admissible $\ZZ/2\ZZ$-covers. We explicitly determine the resulting dual complex, and we use this to define a graph complex which computes the weight zero compactly supported cohomology of $\cH_{g, n}$. Using this graph complex, we give a sum-over-graphs formula for the $S_n$-equivariant weight zero compactly supported Euler characteristic of $\cH_{g, n}$. This formula allows for the computer-aided calculation, for each $g\le 7$, of the generating function $\mathsf{h}_g$ for these equivariant Euler characteristics for all $n$. More generally, we determine the dual complex of the boundary in any moduli space of pointed admissible $G$-covers of genus zero curves, when $G$ is a cyclic group of prime order, as a symmetric $\Delta$-complex. We use these complexes to generalize our formula for $\mathsf{h}_g$ to moduli spaces of $n$-pointed smooth prime cyclic covers of genus zero curves. \textit{This version contains corrections to the published version: see the attached erratum.}
    
    % We also prove that the weight zero compactly supported cohomology $W_0 H^i_c(\cH_{g,n} ; \QQ)$ is supported in degrees $i \geq g$, and write down an explicit cycle spanning $W_0 H^{2g + 1}_c(\cH_{g,2};\QQ)$ for all $g \geq 2$. 
\end{abstract}
\section{Introduction}
For integers $g\ge 2$ and $n\ge 0$, let $\cH_{g,n}\subset \cM_{g,n}$ denote the complex moduli stack of $n$-marked smooth hyperelliptic curves of genus $g$. This space is a smooth Deligne--Mumford stack of dimension $2g + n - 1$. The group $S_n$ acts on $\cH_{g,n}$ by permuting the marked points, and the rational cohomology groups with compact support $H^i_c(\cH_{g,n};\QQ)$ are $S_n$-representations in the category of mixed Hodge structures over $\QQ$. In particular, each cohomology group $H^i_c(\cH_{g,n};\QQ)$ carries a weight filtration
\[ W_0 H^i_c(\cH_{g,n};\QQ) \subset W_1 H^i_c(\cH_{g,n};\QQ) \cdots \subset W_{4g + 2n - 2} H^i_c(\cH_{g,n};\QQ) = H^i_c(\cH_{g,n};\QQ), \]
which is preserved by the $S_n$-action. In this paper, we study the $S_n$-representation defined by the weight zero piece of this filtration.

When $X$ is a smooth and separated variety or Deligne-Mumford stack, Deligne's weight spectral sequence \cite[\S 3.2]{Deligne-hodgeII} computes the associated graded pieces of the weight filtration on the compactly supported cohomology of $X$. It identifies the weight zero piece with the reduced cohomology of the dual complex of any normal crossings compactification of $X$. We will furnish a normal crossings compactification of $\cH_{g, n}$ using the theory of pointed admissible $\ZZ/2\ZZ$-covers, as developed by Abramovich--Vistoli \cite{abramovich-vistoli-stable-maps}, Abramovich--Corti--Vistoli \cite{abramovich-corti-vistoli-twisted}, and Jarvis--Kaufmann--Kimura \cite{jarvis-kaufmann-kimura-pointed}, following Harris--Mumford's original theory \cite{harris-mumford-kodaira}. Denoting the dual complex of the resulting boundary divisor by $\Theta_{g,n}$, we then study the weight zero compactly supported cohomology of $\cH_{g,n}$ via the identification
\begin{equation}\label{eqn:weightzero-identification}
W_0 H^i_c(\cH_{g,n};\QQ) \cong \widetilde{H}^{i - 1}(\Theta_{g,n};\QQ)
\end{equation}
mentioned above, where $\widetilde{H}^*$ denotes reduced cohomology. Along the way, we also explicitly determine the dual complex of the boundary in any space of pointed admissible $G$-covers of genus zero curves, for prime cyclic groups $G$ (Theorem \ref{thm:boundary-complex}). 

Our main result concerns the $S_n$-equivariant weight zero compactly supported Euler characteristic
\[ \chi^{S_n}\!\left( W_0 H^*_c(\cH_{g, n};\QQ) \right) := \sum_{i = 0}^{4g + 2n - 2} (-1)^i \ch_n\left(W_0 H^i_c(\cH_{g,n};\QQ)\right) \in \Lambda, \]
where $\ch_n(\cdot)$ denotes the \textit{Frobenius characteristic} of an $S_n$-representation: this is an element of the ring \[\Lambda = \lim_{\longleftarrow} \QQ[x_1, \ldots, x_n]^{S_n} \] of symmetric functions, which encodes the character of the representation. See \cite{macdonald-symmetricfunctions} or \cite{Stanley-enumerativecombo2} for more on symmetric functions and the Frobenius characteristic.  

% in the following sense: if $V$ is a rational $S_n$-representation with the irreducible decomposition
% \[V \cong \bigoplus_{\lambda \vdash n} W_{\lambda}^{\oplus a_\lambda}, \]
% where $W_\lambda$ denotes the Specht module corresponding to $\lambda$, then
% \[ \ch_n(V) = \sum_{\lambda \vdash n} a_\lambda s_{\lambda}, \]
% where $s_{\lambda} \in \Lambda$ denotes the \textit{Schur function} corresponding to $\lambda$; the Schur functions for $\lambda \vdash n$ form a basis, orthonormal with respect to the Hall inner product, of the homogeneous degree $n$ part of $\Lambda$. See \cite{stanley-introduction} or \cite{macdonald} for more on symmetric functions and the Frobenius charactteristic.  

% \maddie{here we should either say more, point to a future defnition, or give a citation. Something like paragraphs 2-3 of 1904.06367}. 

For each $g \ge 2$, we define
\[ \mathsf{h}_g := \sum_{n \geq 0}\chi^{S_n}\!\left( W_0 H^*_c(\cH_{g, n};\QQ) \right) \]
to be the generating function for these equivariant Euler characteristics. Note that $\mathsf{h}_g$ is an element of $\hat{\Lambda}$, the degree completion of $\Lambda$. In Theorem \ref{mainthm:FrobChar} below, we prove a sum-over-graphs formula for the generating function $\mathsf{h}_g$. The precise definition of the terms in the formula can be found in Section \ref{section: graphsumformula}. For now, we only remark that $T_{2g + 2}^{<3}$ is a finite set of trees, and given such a tree $C$ there is a canonically associated vertex-weighted graph $P_C$ which can roughly be understood as a ``tropical double cover" of $C$: see Section \ref{section:compactifications} for details on this perspective.
\begin{customthm}{A}\label{mainthm:FrobChar}
    We have
    \[ \mathsf{h}_g = \sum_{C \in T_{2g + 2}^{<3}} %\mathsf{h}_{P_C, E_C, \mathrm{Aut}(P_C)}, 
      \frac{(-1)^{|E_C|}}{|\Aut(P_C)|} \sum_{\tau \in \Aut(P_C)} \mathrm{sgn}(\tau|_{E_C}) \prod_{k \geq 1} (1 + p_k)^{f(P_C, \tau, k)}
    \]
    where $E_C$ is the set of edges of the tree $C$, 
    $p_k = \sum_{n > 0} x_n^{k} \in \hat{\Lambda}$ is the $k$th power sum symmetric function, and $k \cdot f(P_C, \tau, k)$ is given by the compactly supported Euler characteristic of the set of points in $P_C$ which have orbit of length $k$, under the action of $\tau$.
\end{customthm}

Implementing Theorem \ref{mainthm:FrobChar} on a computer, we are able to compute $\mathsf{h}_g$ explicitly for $2 \le g \le 7$: see Table \ref{table:frobenius_data}. The code is available at \cite{code}. Our data allows us to extract the polynomials $F_n(t) \in \QQ[t]$, for each $n\le 9$, which have the property that
$F_n(g) = \chi^0_c(\cH_{g,n})$ for each $g\ge 2$, where
\[\chi^0_c(\cH_{g,n}) := \sum_{i = 0}^{4g + 2n - 2} (-1)^i \dim_{\QQ} W_0H^i_c(\cH_{g, n};\QQ) \]
denotes the numerical weight zero compactly supported Euler characteristic. See Proposition \ref{prop:euler_char} in Section~\ref{sec: point-counting} below. Also see Figure \ref{fig:h2exmample} in Appendix \ref{sec:calculations} for an illustration of Theorem \ref{mainthm:FrobChar} when $g = 2$; in this case $T_6^{<3}$ consists of three trees and their contributions to $\mathsf{h}_2$ can be computed by hand.

%\ref{conj:euler_char} about the numerical weight zero compactly supported Euler characteristic of $\mathcal{H}_{g,n}$ for $4 \leq n \leq 9$. 

Our proof of Theorem \ref{mainthm:FrobChar} relies on our description of the cellular chain complex of $\Theta_{g, n}$ as a graph complex generated by certain double covers of trees, which are a special case of the theory of graph-theoretic admissible covers we develop in Section \ref{section:dual_complex}. We find that several subcomplexes of this graph complex are acyclic; the proofs are given in Section \ref{section:acyclic}.  As in earlier work on $\cM_{g,n}$ \cite{cgp-marked}, one conceptually important subcomplex is the {\em repeated marking subcomplex}, i.e., the subcomplex spanned by graph-theoretic admissible covers containing a vertex supporting more than one marking. This subcomplex is acyclic (Theorem~\ref{RepAcyclic}), and after quotienting by it, the resulting chain complex is related to {\em configuration spaces of distinct points} on graph-theoretic admissible covers; see~\cite{bibby-chan-gadish-yun-homology-to-appear, bibby-chan-gadish-yun-serre} for related work. Since Theorem~\ref{mainthm:FrobChar} is about Euler characteristics, we may work one graph-theoretic admissible cover at a time, summing the individual contributions.  For each individual graph-theoretic admissible cover, we use Proposition~\ref{CWContribution}, explained more below, to calculate its contribution. This proves Theorem~\ref{mainthm:FrobChar}.  

Proposition~\ref{CWContribution} may be useful in other applications, so we mention it briefly here: it gives a formula for the completed symmetric function
\[\sum_{n \geq 0} \chi_c^{S_n}\!
\left( \left(\mathrm{Conf}_n(X) \times \Delta^\circ \right)\!/G\right),\]
where $X$ is any finite CW complex, $\Delta^\circ$ is an open simplex, $G$ is a finite group, and $G$ acts on $X$ cellularly and on $\Delta^\circ$ by permuting vertices.  See Section~\ref{section: graphsumformula}. This proposition is closely inspired by a result of Gorsky~\cite[Theorem 2.5]{gorsky-equivariant} concerning complex quasi-projective varieties $X$ with an action of a finite group; our specific formulation is a new contribution. In particular, it does not appear in the work of Chan--Faber--Galatius--Payne on the top weight cohomology of $\cM_{g,n}$, where an alternate argument, which is less geometric, is used \cite[Proposition 3.2]{cfgp-sn}.

Now let us turn our attention to individual cohomology groups, rather than Euler characteristics.  First, for $n=0,1,2,$ and $3$, the cohomology of $\cH_{g,n}$ was completely computed by Tommasi \cite{tommasi-thesis}; see  Section~\ref{sec:related}.
The consequences of these computations for the weight zero part of cohomology with compact supports can be interpreted via our work as statements about chain complexes of graph-theoretic admissible covers.  In Section~\ref{section:acyclic}, we prove some of these statements, using the acyclicity results mentioned above. In particular, we deduce the following facts, first proved by Tommasi:
\begin{customprop}{B}\label{mainthm:smalln}
 For all $g \geq 2$, we have 
\begin{enumerate}
\item $W_0 H^i_c(\cH_{g, n}; \QQ) = 0$ for all $i$, when $n \leq 1$;
\item When $n = 2$, we have 
\[ W_0 H_{c}^{2g + 1}(\cH_{g, 2} ;\QQ)  \cong \QQ. \]
As an $S_2$-representation, we have
\[W_0 H_{c}^{2g + 1}(\cH_{g, 2} ;\QQ)  \cong \begin{cases}
\mathrm{triv} &\text{ if } g \text{ is even} \\
\mathrm{sgn} &\text{ if }g \text{ is odd}.
\end{cases}
\]
\end{enumerate}
\end{customprop}
%We prove Proposition~\ref{mainthm:smalln} in Section \ref{section:acyclic}. 
Part (1) of Proposition~\ref{mainthm:smalln} is established via a spectral sequence argument, similar to the ones we use for acyclicity of other subcomplexes of $\Theta_{g, n}$. For part (2), we write down an explicit cellular cycle on $\Theta_{g, 2}$ corresponding to the nonzero class in $W_0 H^{2g + 1}(\cH_{g, 2};\QQ)$: see Figure \ref{fig:H_g2-cycle} in Section \ref{section:acyclic}. Tommasi shows additionally that $W_0 H^i_c(\cH_{g, 2};\QQ) = 0$ for $i \neq 2g + 1$, but we do not see how to prove this directly using our graph complex, nor have we investigated whether we can use our methods to re-deduce $W_0 H^*_c(\cH_{g,3};\QQ)$ for all $g$.

\subsection{The support of $W_0H^*_c(\cH_{g, n};\QQ)$} 
It is worth noting that the weight zero compactly supported cohomology of $\cH_{g,n}$ is supported in at most two degrees. Precisely, \begin{equation}\label{eq:the-two}
    W_0 H^i_c(\cH_{g,n};\QQ)=0 \quad \text{ unless } \quad i=2g-2+n \text{ or }i=2g-1+n.
\end{equation} %is supported in only two degrees $i=2g-2+n$ and $i=2g-1+n.$   
 %Indeed, as explained by D.~Petersen in a MathOverflow post \cite{petersen-overflow-Hgn}, the \textit{affine stratification number} \cite{roth-vakil} of $H_{g, n}$ is $1$ for all $n > 0$: this statement is true for $n=1$ because $H_{g,1}$ is stratified by the two affine loci in which the marked point is, respectively is not, a Weierstrass point.  Then statement then follows for all $n\ge 1$ because the morphism $\cH_{g,n+1}\to \cH_{g,n}$ is affine for $n\ge 1$.  
We now explain the claim~\eqref{eq:the-two}, which follows from an argument we learned from D.~Petersen.  To sidestep stack-theoretic issues, let us momentarily replace $\cH_{g, n}$ by its coarse moduli space $H_{g, n}$; this is inconsequential on the level of rational cohomology. It is well-known that $H_{g}$ is affine, as it can be identified with the quotient $\cM_{0, 2g+2}/S_{2g + 2}$. In general, $H_{g, n}$ is not far from affine: as explained by D. Petersen in a MathOverflow post \cite{petersen-overflow-Hgn}, the \textit{affine stratification number} \cite{roth-vakil} of $H_{g, n}$ is $1$ for all $n > 0$. By \cite[Corollary 4.19]{roth-vakil} and a suitable comparison theorem for {\'etale cohomology} \cite[Theorem 21.1]{milneLEC}, we may conclude that
\[H^i(\cH_{g, n}; \QQ) = 0 \text{ for }i> 2g + n,  \quad \text{and} \quad 
%and by Poincar\'e duality, 
H^i_c(\cH_{g, n}; \QQ) = 0 \text{ for }i < 2g -2 + n,  \]
the latter by Poincar\'e duality.
As the dual complex $\Theta_{g, n}$ of the normal crossings compactification of $\cH_{g, n}$ by pointed admissible $\ZZ/2\ZZ$-covers is a generalized cell complex of dimension $2g - 2 + n$ (Section~\ref{section:dual_complex}), the claim~(\ref{eq:the-two}) follows immediately from (\ref{eqn:weightzero-identification}).
%we may conclude that
%\[ W_0H^{i}_c(\cH_{g, n};\QQ) = 0 \]
%for $i \notin \{2g - 2 + n, 2g - 1 + n\}$, by (\ref{eqn:weightzero-identification}). 
%In particular, it follows that
%\[ \chi^{S_n} \left( W_0 H^*_c(\cH_{g, n};\QQ)\right) = (-1)^n \left(\ch_n(W_0 H^{2g -2 + n}_c(\cH_{g, n};\QQ)) - \ch_n(W_0 H^{2g -1 + n}_c(\cH_{g, n};\QQ))\right).  \]

Thus, our formula for $\mathsf{h}_g$ is a formula for the {\em difference} of the two $S_n$-representations in~\eqref{eq:the-two} and can be used to bound the multiplicities of Specht modules appearing in them individually.  We have not investigated whether  $\mathsf{h}_g$ is in fact a cancellation-free formula for this difference.

% We recall his argument here. One first stratifies \[H_{g, 1} = H_{g, 1}^{u} \coprod H_{g, 1}^w,\] where $H_{g, 1}^{w}$ is the locus where the marked point is equal to a Weierstrass point, and $H_{g, 1}^u$ is the complement of $H_{g, 1}^w$. We have identifications
% \[ H_{g, 1}^w \cong \cM_{0, 2g + 2}/S_{2g + 1} \quad\mbox{and}\quad H_{g, 1}^u \cong \cM_{0, 2g + 3}/S_{2g + 2}, \]
% so both loci are affine. Since $H_{g, 1}$ is not itself affine (any fiber of the map $H_{g, 1}\to H_g$ is a complete curve in $H_{g, 1}$), one concludes that the affine stratification number of $H_{g, 1}$ is $1$. Finally, one observes that the morphism $H_{g, n} \to H_{g, 1}$ which erases all but the first marked point is an affine morphism, so the affine stratification number of $H_{g, n}$ is also $1$ for all $n$. 

\subsection{Related work on the cohomology of $\cH_{g, n}$}\label{sec:related}
Recently, there have been a number of significant advances on the geometry of moduli spaces of pointed hyperelliptic curves.  Canning--Larson study the rational Chow ring of $\cH_{g,n}$, in particular determining it completely for $n\le 2g+6$ \cite{canning-larson-rational}.  Their results also have implications for rationality of $\cH_{g,n}$. More generally, there has been progress on understanding the birational geometry of $\cH_{g,n}$; see, for example, the overview and references in that paper.  In another direction, Bergstr\"om--Diaconu--Petersen--Westerland \cite{bergstrom-diaconu-petersen-westerland-hyperelliptic} compute the stable homology of braid groups with coefficients in (any Schur functor applied to) the Burau representation. These results have implications for the stable homology of moduli spaces of hyperelliptic curves with twisted coefficients. They can also be related to the Serre spectral sequence on rational cohomology for the fiber bundle $\Conf_n(S_g) \to \cH_{g,n} \to \cH_g$, as C.~Westerland has explained to us. 
Our focus here is the cohomology groups of $\cH_{g,n}$ with (untwisted) $\QQ$-coefficients, and specifically the weight zero compactly supported cohomology groups.

% For a small number of marked points, Tommasi completely determines the weight filtration on rational cohomology in her thesis  \cite{tommasi-thesis}. She shows that $\cH_{g}$ has the rational cohomology of a point, and determines the weight filtration on the cohomology of $\cH_{g, 2}$, with its $S_2$-action. D. Peterson has communicated to us that her methods also apply to determine the weight filtration on the cohomology of $\cH_{g, n}$ for $n =1$ and $3$.

The topological Euler characteristic of $\cH_{g, n}$ has been computed by Bini \cite{bini-Hgn}, but his techniques are not compatible with the weight filtration. Gorsky \cite{gorsky-hyperelliptic} calculates the equivariant Euler characteristic
\[ \chi^{S_n}(\cH_{g, n}) := \sum_{i = 0}^{4g + 2n - 2} (-1)^i\ch_n(H^{i}(\cH_{g, n};\QQ)), \]
by fibering $\cH_{g, n}$ over $\cH_{g}$. The fiber of this morphism over a point of $\cH_{g}$ representing a curve $C$ is equal to $\mathrm{Conf}_n(C)/\!\Aut(C)$. Gorsky proceeds by stratifying $\cH_{g}$ by the $S_n$-equivariant Euler characteristic of the fibers, and then calculating the Euler characteristic of each stratum. Our techniques are similar in spirit to Gorsky's. The $S_n$-equivariant weight zero compactly supported Euler characteristic of $\cH_{g, n}$ is equal to $h_n - \chi^{S_n}(\Theta_{g, n})$, where $h_n \in \Lambda$ is the $n$th homogeneous symmetric function. As explained above, we first remove an acyclic locus from $\Theta_{g, n}$, and then stratify the remaining space in terms of configuration spaces of graphs, summing up these contributions to give our formula (Section \ref{section: graphsumformula}).

\subsection{Relation to point-counting}\label{sec: point-counting} Bergstr\"{o}m \cite{bergstrom-equivariant} studies the cohomology of $\cH_{g, n}$ via point-counting: for all $g \ge 2$, he gives an algorithm to determine the count of $\FF_q$-points of $\cH_{g, n}$ for $n \leq 7$ and for all prime powers $q$. Together with the results of \cite{BEGR-pointcounts}, Bergstr\"{o}m's work implies that for odd $q$, the number of $\FF_q$-points of $\cH_{g, n}$ agrees with a polynomial $P_{g, n}(q)$ for $n \leq 9$ (there is a different polynomial for even $q$). By \cite[Theorem 6.1.2(3)]{hrk-point-counting}, we have an equality
\[P_{g, n}(q) = \sum_{j = 0}^{2g + n - 1} \chi_c^{2j}(\cH_{g,n}) q^j, \]
where
\[ \chi_c^k(\cH_{g,n}) :=  \sum_{i = 0}^{4g + 2n - 2} (-1)^i\dim_{\QQ} \mathrm{Gr}_k^{W} H^i_c(\cH_{g, n};\QQ), \]
and
\[\mathrm{Gr}_k^{W} H^i_c(\cH_{g, n};\QQ) := W_k H^i_c(\cH_{g, n};\QQ) / W_{k - 1} H^i_c(\cH_{g, n};\QQ) \]
is the $k$th associated graded piece of the weight filtration. In particular, the constant term of $P_{g, n}(q)$ is equal to the weight zero compactly supported Euler characteristic. Bergstr\"{o}m's original work \cite{bergstrom-equivariant} is $S_n$-equivariant, and we have confirmed that our data agrees with his for $n \leq 7$. He has explained to us that \cite[Theorem 5.2]{bergstrom-equivariant} and \cite{BEGR-pointcounts} imply that for each $n \leq 9$, there exists a polynomial $F_n(t) \in \QQ[t]$, with degree bounded by $n - 2$ if $n$ is even and $n - 3$ if $n$ is odd, such that
\[ \chi_c^0(\cH_{g, n}) = F_n(g) \]
for all $g$. With these bounds on the degrees, our formula allows us to compute this polynomial for all $n \leq 9$, using the data in Table \ref{table:euler_char}. The polynomials $F_n(t)$ can certainly be calculated from Bergstr\"om's work, but did not explicitly appear there, so we record them below. In each case, the degree of $F_n(t)$ attains the communicated bound. 
% \newpage
{\begin{customprop}{C}
    \label{prop:euler_char}
    We have $\chi^0_c(\cH_{g, n}) = 0$ for $n \in \{0, 1, 3\}$, while $\chi_c^0(\cH_{g, 2}) = -1$. For $4 \leq n \leq 9$, we have the following:
    
    {\renewcommand{\arraystretch}
    {1.5}
    \centering
    \begin{tabular}{l}
$\chi^0_c(\mathcal{H}_{g,4}) = g(1-g)$                                                                 \\
$\chi^0_c(\mathcal{H}_{g,5}) =5g(-1+g)$                                                                \\
$\chi^0_c(\mathcal{H}_{g,6}) = \frac{1}{8} g (198 - 203 g + 18 g^2 - 13 g^3)$                          \\
$\chi^0_c(\mathcal{H}_{g,7}) =\frac{7}{4} g (-78 + 83 g - 18 g^2 + 13 g^3)$                            \\
$\chi^0_c(\mathcal{H}_{g,8}) = \frac{1}{4} g (3420 -3784 g + 1355 g^2 - 1005 g^3 + 25 g^4 - 11 g^5)$   \\
$\chi^0_c(\mathcal{H}_{g,9}) = \frac{9}{4} g (-2700 + 3092 g - 1545 g^2 + 1195 g^3 - 75 g^4 + 33 g^5)$.
\end{tabular}\par}
    
\end{customprop}}
% We conjecture that the polynomials in $g$ computing $\chi_c^0(\cH_{g,n})$ for each fixed $n$ have vanishing constant term, and that the coefficients are alternating and have unimodal absolute values.
% siddarth: is it even clear that that there exists a polynomial for n>9? leaving this out for now
% 
% and are of degree $\lfloor n/2\rfloor + 1$
%  siddarth: it seems like jonas already knows that a bound on the degree is n-2 or n-3 depending on the parity of $n$, and this bound is obtained in all of the examples above, so this is equally natural to conjecture imo

\subsection{Relation to previous work on $\cM_{g, n}$} Our calculations are a new step in understanding weight zero compactly supported rational cohomology of moduli spaces via combinatorics of normal crossings compactifications %$\cM_g$ and $\cM_{g,n}$ 
\cite{acp, bbcmmw-top, cgp-graph-homology, cgp-marked, cfgp-sn}. %, and moduli spaces of abelian varieties $\cA_g$ \cite{bbcmmw-top}. 
In our calculation of $\mathsf{h}_g$, we proceed in a similar fashion to Chan--Faber--Galatius--Payne \cite{cfgp-sn}, who calculate the $S_n$-equivariant weight zero Euler characteristic of $\cM_{g, n}$. They use the dual complex $\Delta_{g, n}$ of the Deligne--Mumford--Knudsen compactification $\cM_{g, n} \subset \overline{\cM}_{g, n}$, which can be interpreted as a tropical moduli space of curves \cite{acp}. They express the generating function
\[\mathsf{z}_g := \sum_{n \geq 0} \chi^{S_n}\!\left( W_0 H^*_c(\cM_{g, n}; \QQ)\right) \]
as a sum over contributions from configuration spaces of graphs. The contribution from each graph is a sum of monomials in the inhomogeneous power sum symmetric functions $P_i := 1+p_i$, of degree equal to the topological Euler characteristic of the graph. A crucial difference between their work and ours, which has been an unexpected subtlety here, is that they find that the only graphs contributing to their formula are connected with first Betti number $g$. As such, their formula for $\mathsf{z}_g$ is a Laurent polynomial in the $P_i$'s, homogeneous of degree $1 - g$, where $P_i$ has degree $i$. The ability to focus on graphs with fixed Euler characteristic is a significant conceptual aid to their work. In contrast, we find that while all of the graphs contributing to $\mathsf{h}_g$ are connected double covers of metric trees, they do not have fixed first Betti number, so their topological Euler characteristics vary, and indeed for $g \geq 3$ the formulas for $\mathsf{h}_g$ are not homogeneous in the $P_i$'s. When $g = 2$, we have $\cH_{2, n} = \cM_{2,n}$, so $\mathsf{h}_2 = \mathsf{z}_2$ is homogeneous of degree $-1$.

\subsection{Applications to moduli spaces of admissible $G$-covers in genus zero}
While our main focus in this paper is the moduli space $\cH_{g, n}$, our techniques are more general. As mentioned above, Theorem \ref{thm:boundary-complex} in Section \ref{section:dual_complex} contains a description of the dual complex of the boundary divisor in any moduli space of pointed admissible $G$-covers of genus zero curves, when $G$ is a prime cyclic group. We specialize to $G = \ZZ/2\ZZ$ in order to study $\cH_{g, n}$. We can prove a generalization of Theorem \ref{mainthm:FrobChar} to more general moduli spaces of pointed $G$-covers: see Remarks \ref{remark: generalizedrepeatedmarkings} and \ref{remark: generalizedgraphsum}, and Theorem \ref{mainthm:extendedFrobChar} in Section \ref{section: graphsumformula}.

\subsection*{Acknowledgments} We are grateful to Dan Abramovich for teaching us about twisted stable maps and admissible $G$-covers, and to Jonas Bergstr\"{o}m for explaining his work \cite{bergstrom-equivariant} and sharing his data on the weight zero compactly supported Euler characteristic of $\cH_{g, n}$. Jonas Bergstr\"{o}m, Dan Petersen, and Dhruv Ranganathan provided extremely valuable comments on a draft of this paper; we thank them very much.  Finally, we sincerely thank two anonymous referees who combed through our manuscript and provided incisive and thorough comments throughout.
MB is supported by the National Science Foundation under Award No. 2001739. MC was supported by NSF CAREER DMS-1844768, a Sloan Foundation Fellowship and a Simons Foundation Fellowship. SK was supported by an NSF Graduate Research Fellowship. Any opinions, findings, and conclusions or recommendations expressed in this material are those of the authors and do not necessarily reflect the views of the National Science Foundation.

\section{Pointed admissible $G$-covers and their moduli}

In this section we recall moduli spaces of pointed admissible $G$-covers, following \cite{abramovich-vistoli-stable-maps, abramovich-corti-vistoli-twisted,harris-mumford-kodaira, jarvis-kaufmann-kimura-pointed}.
We determine the connected components of these spaces when $g=0$ and $G$ is prime cyclic (Proposition \ref{prop:cc-smooth-case}), and give a normal crossings compactification (Proposition \ref{prop:normal_crossings}). Later, in Section \ref{section:dual_complex}, we will determine the dual complex of this compactification. Ultimately, we will obtain a normal crossings compactification of $\mathcal{H}_{g,n}$ and the corresponding dual complex as a special case
in Section~\ref{section:compactifications}.

\subsection{Admissible $G$-covers}
Let $G$ be a finite group, and let $g,n\ge 0$ be integers such that ${2g-2+n>0}$.
We recall the notion of an {\em admissible $G$-cover} of nodal curves of type $(g,n)$ over an arbitrary base scheme $T$ (\cite[Definition 2.1]{jarvis-kaufmann-kimura-pointed}, \cite[Definition 4.3.1]{abramovich-corti-vistoli-twisted}). %\cite[Definition 2.1]{jarvis-kaufmann-kimura-pointed}).  
It is the data of an $n$-marked, stable genus $g$ curve $(C,p_1,\ldots,p_n)$ over $T$, and a covering of nodal curves $\phi\col P\to C$ with an action of $G$ on $P$ leaving $\phi$ invariant, such that:
\begin{enumerate}
    \item $\phi$ is a principal $G$-bundle away from the nodes and markings of $C$,
    \item The analytic local equations for $P\to C\to T$ at a point $p\in P$ over a node of $C$ are
    \[\Spec A[z,w]/(zw-t) \to \Spec A[x,y]/(xy-t^r)\to \Spec A,\]
    where $t\in A$, $x = z^r$ and $y=w^r$ for some integer $r>0$.
    \item The analytic local equations for $P\to C\to T$ at a point $p\in P$ over a marked point of $C$ are
    \[\Spec A[z]\to \Spec A[x] \to \Spec A,\]
    where $x=z^s$ for some integer $s>0$.
    \item if $x\in P$ is a geometric node, then the action of the stabilizer $G_x$ of $x$ on the tangent spaces of the two analytic branches at $x$ is {\em balanced}: the characters of these two one-dimensional representations of $G_x$ are inverse to each other.
\end{enumerate}

\noindent 
Admissible $G$-covers of type $(g,n)$ form a Deligne-Mumford stack, denoted $\Adm_{g,n}(G)$; this is a consequence of the  identification of $\Adm_{g,n}(G)$ with the space $\cB^{\mathrm{bal}}_{g,n}(G)$ of {\em balanced twisted $G$-covers} of type $(g,n)$ %in \cite[Theorem 4.3.2]{abramovich-corti-vistoli-twisted} with the space $\cB_{g,n}^{\mathrm{bal}}$, 
which is proven in \cite{abramovich-vistoli-stable-maps} to be a Deligne-Mumford stack.  We may write {\em $G$-cover} rather than {\em admissible $G$-cover} for short.
%\cite[Theorem 4.3.2]{abramovich-corti-vistoli-twisted}.

\subsection{Admissible covers of smooth curves}

Let $\Adm_{g,n}^\circ(G)$ denote the open substack of %admissible 
$G$-covers in which the target curve (and hence also the source curve) is smooth. In this section, we will determine the connected components of $\Adm_{0,n}^\circ(G)$ (Proposition~\ref{prop:adm-m0n-trivial}). We will use this result later when determining the connected components of the corresponding space of pointed admissible $G$-covers.

%; the objects parametrized may simply be called {\em $G$-covers} of smooth curves in lieu of the longer ``admissible $G$-covers.''
 There is a forgetful map \[\pi \col \Adm_{g,n}^\circ(G)\to \cM_{g,n}\] sending a $G$-cover $P\to (C,p_1,\ldots,p_n)$  to the $n$-pointed curve $(C,p_1,\ldots,p_n)$.  The morphism $\pi$ is \'etale; this property can be deduced from %the description of the map of formal deformation spaces induced by $\pi$, in 
 \cite[Theorem 5.1.5]{bertin-romagny-champs}, as explained in Proposition 6.5.2 of op.~cit.  Working over $\CC$, the fiber over  
 $(C,p_1,\ldots,p_n)$ %$\pi^{-1}(C,p_1,\ldots,p_n)$ 
 is identified with the set 
 \begin{equation}\label{eq:set}\Hom(\pi_1(C-\{p_1,\ldots,p_n\},p_0), G)/G\end{equation}
 where $G$ acts by conjugation, and $p_0 \in C-\{p_1,\ldots,p_n\}$ is any choice of basepoint. 
 There are no other restrictions on the set~\eqref{eq:set}; in particular, the source curves $P$ are not required to be connected.
 An element of the set~\eqref{eq:set} specifies a $G$-cover of the punctured curve $C-\{p_1,\ldots,p_n\}$, which can be extended uniquely over the punctures.  Then the data of the morphism $\pi$ is equivalent to the data of the action of $\pi_1$ of the base space $\cM_{g,n}$ on the fiber~\eqref{eq:set} above. 
 We shall now consider this action in the case $g=0$, when the action may be understood via the classical Hurwitz theory  of $\PP^1$.  We denote by \[\varepsilon^\mathrm{ni}_n(G) := \{(g_1,\ldots,g_n)\in G^n: g_1\cdots g_n = 1\}\]
the set of {\em Nielsen classes}.  We do not impose that $g_1,\ldots,g_n$ generate $G$;  correspondingly, our source curves are not required to be connected. The group $G$ acts by conjugation on $\varepsilon^\mathrm{ni}_n(G)$, and the elements of $\varepsilon^\mathrm{ni}_n(G)/G$ are called {\em inner Nielsen classes}.  Recall the following relationship between the set~\eqref{eq:set} to the set of inner Nielsen classes: choose loops  $\rho_1,\ldots,\rho_n$ around $p_1,\ldots,p_n$, respectively, based at $p_0$, such that $\rho_1,\ldots,\rho_n$ generate $\pi_1(C-\{p_1,\ldots,p_n\},p_0)$ subject only to the relation \[\rho_1\cdot\cdots\cdot \rho_n = 1.\] Such a choice identifies the set~\eqref{eq:set} with the inner Nielsen classes.

Now the following diagram of pullback squares relates $\Adm_{0,n}^\circ(G)$ to Hurwitz spaces of $G$-covers.
\[\xymatrix@R=6mm@C=6mm{
\Adm_{0,n}^\circ(G) \ar[d]\ar[r] & \cH_{\PP^1,n}^G \ar[d] \ar[r] & U\cH^G_{\PP^1,n} \ar[d] \\ \cM_{0,n} \ar[r] &\Conf_n(\PP^1)\ar[r] & \mathrm{UConf}_n(\PP^1)
}\]
The spaces above are defined as follows. The configuration spaces (ordered and unordered) of $n$ points in $\PP^1$ are denoted $\Conf_n(\PP^1)$ and $\mathrm{UConf}_n(\PP^1)$, respectively.
The space $U\cH^G_{\PP^1,n}$ is the moduli space parametrizing sets $S\subset \PP^1$ of $n$ points, together with a ramified $G$-cover $f\col P\to \PP^1$ whose branch locus is contained in $S$.  The space $\cH^G_{\PP^1,n}$ is the ordered version of this space, obtained by pullback. The map $\cM_{0,n}\to \Conf_n(\PP^1)$ fixes $(p_1,p_2,p_3)$ to be $(0,1,\infty)$, for instance.

 \propnow{\label{prop:adm-m0n-trivial}If $G$ is abelian, then $\cH^G_{\PP^1,n}\to \mathrm{Conf}_n(\PP^1)$, and hence also $\Adm_{0,n}^\circ(G)\to \cM_{0,n}$, is a trivial bundle.  As a variety, $\Adm_{0,n}^\circ(G)$ is isomorphic to $\cM_{0,n} \times \varepsilon^{\mathrm{ni}}_n(G)$.} 

 \begin{proof}
For an arbitrary finite group $G$, the way in which $U\cH_{\PP^1,n}^G$ is a covering space over $\mathrm{UConf}_n(\PP^1)$ is classically understood, essentially going back to Hurwitz \cite{hurwitz-ueber}, see \cite[p.~547]{fulton-hurwitz}.
The following is a complete description. 
%% added:
Let $S = \{s_1,\ldots,s_n\}$.  
For an appropriate choice of basis, $\Hom(\pi_1(\PP^1-S,p_0),G)$ is identified with $\varepsilon^{\mathrm{ni}}_n(G)$.  And $\pi_1(\mathrm{UConf}_n(\PP^1))$ has a presentation with generators $\gamma_1,\ldots,\gamma_{n-1}$, where $\gamma_i$ interchanges points $i$ and $i+1$.  Furthermore the generators $\gamma_i$ act on $\varepsilon^{\mathrm{ni}}_n(G)$ via
\[\gamma_i \cdot (g_1,\ldots,g_n) = (g_1,\ldots,g_{i-1}, g_i g_{i+1} g_i^{-1}, g_{i}, g_{i+2},\ldots, g_n). \]
% The way in which $\Adm_{0,n}^\circ(G)$ is a covering space of $\cM_{0,n}$ can be similarly understood, at least in principle.
In the case that $G$ is abelian, the action is 
\[\gamma_i \cdot (g_1,\ldots,g_n) = (g_1,\ldots,g_{i-1}, g_{i+1}, g_{i}, g_{i+2},\ldots, g_n). \]
In other words, the action of $\pi_1(\mathrm{UConf}_n(\PP^1))$ on $\varepsilon^{\mathrm{ni}}_n(G)$ factors through $\pi_1(\mathrm{UConf}_n(\PP^1))\to S_n$.
Passing to the {\em ordered} configuration space, we therefore obtain a trivial action of the spherical braid group $\pi_1(\mathrm{Conf}_n(\PP^1))$ on $\varepsilon_n^{\mathrm{ni}}(G)$, proving the claim.
\end{proof} %Thus we have arrived at the following statement for abelian groups.

%The maps $\cM_{0,n} \to \mathrm{Conf}_n(\PP^1) \to \mathrm{UConf}_n(\PP^1)$,
%to the ordered and unordered configuration spaces of $n$ points on $\PP^1$, respectively, induce maps 

%\[\mathrm{Mod}_{0,n}\to \mathrm{PB}_n(S^2) \to B_n(S^2).\]

\begin{remark}
\label{rem:adm0_stack}
Stack-theoretically, we have
 \begin{equation}\label{eq:trivial-bundle} \Adm_{0,n}^\circ(G) \cong \cM_{0,n} \times [\varepsilon^{\mathrm{ni}}(G)/G]\end{equation} if $G$ is abelian,
 where $G$ acts trivially on $\varepsilon^{\mathrm{ni}}(G)$.  
Under this identification, write
\begin{equation}\label{eq:adm-gi}
\Adm_{0,n}^\circ(G;g_1,\ldots,g_n)
\end{equation}
for the connected component of $\Adm_{0,n}^\circ(G)$ corresponding to the Nielsen class $(g_1,\ldots,g_n)$; it is isomorphic to $\cM_{0,n}\times BG$. %It parametrizes those admissible covers $P\to (C,p_1,\ldots,p_n)$ for which the corresponding map $\pi_1(C-\{p_1,\ldots,p_n\},p_0)\to G$ takes any keyhole loop around $p_i$ to $g_i$, for each $i$.
\end{remark}

\subsection{Pointed admissible covers}
We study spaces of pointed admissible covers and determine the connected components of these spaces in Proposition \ref{prop:cc-smooth-case}. This is an important calculation towards the computation of the connected boundary strata in Theorem \ref{thm:boundary-complex}, since the boundary strata of spaces of pointed admissible covers are quotients of products of smaller spaces of pointed admissible covers.

Let $G$ be any group, not necessarily abelian.  Let $\ov{\cM}_{g,n}^G$ denote the space of $n$-marked {\em pointed admissible $G$-covers}
of genus $g$ \cite{jarvis-kaufmann-kimura-pointed}.  It is a moduli space for nodal admissible $G$-covers $P\to (C,p_1,\ldots,p_n)$, together with a choice of a lift $\widetilde{p_i}$ on $P$ of each $p_i$. 
The open substack $\cM_{g,n}^G$ is the moduli space of pointed admissible $G$-covers in which source and target are smooth. Summarizing, we have a Cartesian square
\[\squarediagramlabel{\cM_{g,n}^G}{\ov{\cM}_{g,n}^G}{\Adm^\circ_{g,n}(G)}{\Adm_{g,n}(G)}{\subset
}{\pi}{\pi}{\subset}\]
which lays out the unfortunate lack of parallelism in the notation for these spaces. The notation comes from the literature, however.

\propnow{\label{prop:etale}The morphisms $\cM_{g,n}^G\to \Adm^\circ_{g,n}(G)$ and $\ov{\cM}_{g,n}^G \to \Adm_{g,n}(G)$ are \'etale.}

\noindent For easy reference, we prove Proposition~\ref{prop:etale} below. We note, however, that the argument appears as part of the proof in \cite[Theorem 2.4]{jarvis-kaufmann-kimura-pointed} of the fact that $\ov{\cM}_{g,n}^G$ is a smooth Deligne-Mumford stack, flat, proper, and quasi-finite over $\ov{\cM}_{g,n}$.  

\proofnow{We verify the second statement, which implies the first.  
Recall the construction of $\ov{\cM}_{g,n}^G$, which we summarize following \cite{jarvis-kaufmann-kimura-pointed}.  Let $E\to \mathcal{C}=[E/G]$ denote the universal source curve and stacky target curves, respectively, over $\Adm_{g,n}(G)$, and let $C$ denote the coarse space of $\cC$. For $i=1,\ldots,n$, let $\cS_i \to \cC$ denote the closed substack of $\cC^\mathrm{sm}$ whose image in $C$ is the universal $i^{\mathrm{th}}$ marked point; $\cS_i$ is an \'etale gerbe over $\Adm_{g,n}(G)$.  Let $E_i = E\times_\cC \cS_i$.  We have the following diagram, whose %top two squares are Cartesian: 
top square is Cartesian and where the morphisms known to be \'etale are labeled:

%\[\xymatrix{E_i \ar[r] \ar[d]^{\textcolor{black}{\textrm{\'et}}} & E \ar[r]\ar[d]^{\textcolor{black}{\textrm{\'et}}} & \bu\ar[d]^{\textcolor{black}{\textrm{\'et}}} \\ \cS_i \ar[r]\ar[rdd]^{\textcolor{black}{\textrm{\'et}}} & \cC \!=\! [E/G] \ar[r]\ar[d] & BG \ar[d] \\ & C \ar[r]\ar[d] & \bu \\ & \Adm_{g,n}(G) & 
%}\]
\[\xymatrix{E_i \ar[r] \ar[d]^{\textcolor{black}{\textrm{\'et}}} & E \ar[d]^{\textcolor{black}{\textrm{\'et}}}  \\ \cS_i \ar[r]\ar[rdd]^{\textcolor{black}{\textrm{\'et}}} & \cC \!=\! [E/G] \ar[d] \\ & C \ar[d] \\ & \Adm_{g,n}(G) 
}\]
The morphism $E_i\to \Adm_{g,n}(G)$ is \'etale since it is a composition of $E_i\to \cS_i$, which is a pullback of an \'etale morphism and hence \'etale, and the \'etale gerbe $\cS_i\to \Adm_{g,n}(G)$.    %Now $\ov{\cM}_{g,n}^G$ is constructed as
Therefore
\[\ov{\cM}_{g,n}^G = E_1\times_{\Adm_{g,n}(G)} \cdots \times_{\Adm_{g,n}(G)} E_n\]
is also \'etale over $\Adm_{g,n}(G)$.
}

%\subsection{Connected components of spaces of admissible covers}

The spaces $\cM_{g,n}^G$ and $\Adm^\circ_{g,n}(G)$ %, and their compactifications $\ov{\cM}_{g,n}^G$ and $\Adm_{g,n}(G)$, 
need not be connected, as observed in Remark \ref{rem:adm0_stack}. Given $g_1,\ldots,g_n\in G$, write
$\cM_{g,n}^G(g_1,\ldots,g_n)$ 
%There is a forgetful map \begin{equation}\label{eq:forget-points}\cM_{g,n}^G \to \mathrm{Adm}_{g,n}^{\circ}(G).\end{equation}
%Recall that \[\cM_{g,n}^G(g_1,\ldots,g_n)\] denotes 
for the open and closed substack of $\cM_{g,n}^G$ in which the monodromy at the marking $\widetilde{p_i}$ in the source curve is $g_i$.  We recall the notion of monodromy at a point in the source curve, following \cite[\S2.1]{jarvis-kaufmann-kimura-pointed}: pick a small oriented loop around the point $p_i$ in the target curve, say based at a point $q_i$ near $p_i$; then the loop lifts to $d$ possible paths between the $d$ preimages of $q_i$ near $\wt p_i$, where $d$ temporarily denotes the number of sheets of $P$ meeting at $\wt p_i$.  Each of these $d$ paths starts and ends at points $x$ and $gx$, respectively, for some well-defined $g\in G$. Indeed, this $g$ is independent of choice of one of those $d$ paths, since they are each of the form $g^ix$ to $g^{i+1}x$ for $i=0,\ldots,d-1$. The monodromy at $\wt p_i$ is then defined to be $g$.  (Note that $g$ {\em can} depend, a priori, on choice of lift $\wt p_i$ of $p_i$, if $G$ is not abelian. Indeed, the action of any $h\in G$ moves the previously mentioned path near $\wt p_i$ from $x$ to $gx$ to a path near $h\wt p_i$ from $hx$ to $hgx = (hgh^{-1})hx$, so the monodromy at $h \wt p_i$ is $hgh^{-1}$.)
% Write $\mathrm{Adm}_{g,n}^{\circ}(G;g_1,\ldots,g_n)$ for the open and closed substack of $\mathrm{Adm}_{g,n}^{\circ}(G)$ in which the monodromy around the marked points $p_i$, relative to an arbitrary choice of basepoint $\widetilde{p_0}$ in $P$, is equal to $(g_1,\ldots,g_n)$ up to conjugation by $G$.

%\lemnow{Let $G$ be an abelian group.  The map $\cM_{g,n}^G \to \mathrm{Adm}_{g,n}^{\circ}(G)$ restricts to a map $\cM_{g,n}^G(g_1,\ldots,g_n)\to\mathrm{Adm}_{g,n}^{\circ}(G;g_1,\ldots,g_n)$.}

%\proofnow{If an admissible $G$-cover $(P\to C, p_1,\ldots,p_n)$ lies in $\cM_{g,n}^G(g_1,\ldots,g_n)$, then it sends a keyhole loop starting and ending at a chosen basepoint $p_0$ and going around $p_i$ to $g_i$.  Then the monodromy around any lift $\widetilde{p_i}$ of $p_i$ is $g_i$ or a conjugate thereof, so is exactly $g_i$ provided $G$ is abelian.}

\propnow{\label{prop:cc-smooth-case} Let $G$ be a prime cyclic group. %right?
Suppose $g_1\cdots g_n = 1$, so that $\cM_{0,n}^G(g_1,\ldots,g_n)$ is nonempty.  The connected components of $\cM_{0,n}^G(g_1,\ldots,g_n)$ are in bijection with orbits of functions
\[\{1,\ldots,n\}\to G/\langle g_1,\ldots, g_n\rangle \]
under left $G$-translation.  }

\begin{proof}
The restriction of the map $\cM_{0,n}^G\xra{\pi} \Adm^\circ_{0,n}(G)$ to $\cM^G_{0,n}(g_1,\ldots,g_n)$ becomes a surjection  
\[\cM^G_{0,n}(g_1,\ldots,g_n)\xra{\pi} \Adm^\circ_{0,n}(G;g_1,\ldots,g_n) \cong \cM_{0,n}\times BG,\] where the last isomorphism was established in Proposition \ref{prop:adm-m0n-trivial}. This morphism is \'etale by Proposition~\ref{prop:etale}.  
%Proposition~\ref{prop:etale} shows that $\cM_{0,n}^G(g_1,\ldots,g_n)$ is an \'etale cover of $\mathrm{Adm}^{\circ}_{0,n}(g_1,\ldots,g_n)$, which is isomorphic to $[\cM_{0,n}/G] = \cM_{0,n}\times BG$ by Proposition~\ref{prop:adm-m0n-trivial}.
%Let $H = \langle g_1,\ldots, g_n\rangle. $  
Now let $P\to (C, p_1,\ldots,p_n)$ be any unpointed admissible cover; the fiber of $\pi$ over it is the action groupoid on all lifts $\wt{p}_1,\ldots, \wt{p}_n$ of  $p_1,\ldots,p_n$ respectively, with the group $G$ acting by simultaneous translation of the $\wt{p}_i$.  The connected components of $\cM_{0,n}^G(g_1,\ldots,g_n)$ are in bijection with the orbits of this category under the further action of pure mapping class group $\mathrm{Mod}_{0,n}$.  Those orbits are in bijection with orbits of functions $\{1,\ldots,n\}\to \pi_0(P)$ under left $G$-translation; %in other words, one may move the points $p_i$ in $C$ to go from any choices of $\wt{p}_i$ to any other choices $\wt_{p}'_i$, if and only if $\wt p_i$ and $\wt p'_i$ are in the same connected component of $P$.  
 and $\pi_0(P)\cong G/\langle g_1,\ldots,g_n\rangle$.
%Now, two lifts $(\wt p_1,\ldots, \wt p_n)$ and $ (\wt p_1',\ldots, \wt p_n')$  of $(p_1,\ldots,p_n)$ are in the same orbit of the action of $\pi_1(\Adm^{\circ}_{0,n}(g_1,\ldots,g_n),P\to C)$ if and only if $\wt p_i$ and $\wt p'_i$ are in the same connected component of $P$ for each $i$. Therefore, we may associate, to any connected component of  $\cM^G_{0,n}(g_1,\ldots,g_n)$, a well-defined function \[\{1,\ldots,n\} \to \pi_0(P) \cong G/\langle g_1,\ldots,g_n\rangle.\]   On the other hand, two lifts of $(p_1,\ldots,p_n)$ give isomorphic pointed admissible covers if and only if they determine functions $\{1,\ldots,n\}\to G/\langle g_1,\ldots, g_n\rangle$ that are in the same $G$-orbit, yielding the claimed bijection. 
\end{proof}

\begin{remark}\label{rem:erratum}
    In the published version of this article, Proposition \ref{prop:cc-smooth-case} is claimed at the generality of arbitrary abelian groups. However, the statement of the proposition is not correct at that level of generality. The proof of Proposition \ref{prop:cc-smooth-case} is correct for prime cyclic groups, though the attached erratum gives an alternative and more detailed argument. In particular, all of the results in this article which are stated for prime cyclic groups also hold for those finite abelian groups for which the statement of Proposition \ref{prop:cc-smooth-case} holds.
\end{remark}

It will be convenient to work with pointed curves labelled by arbitrary finite sets.  Thus let $G$ be a finite group, $S$ a finite set, and $\rho\col S\to G$ any function.  For $g\ge 0$ with $2g-2+|S|>0$,  let
\[\ov{\cM}^G_{g,S}(\rho)\]
denote the space of pointed admissible $G$-covers of genus $g$ curves with specified monodromy $\rho$. Let ${\cM}^G_{g,S}(\rho)$ denote the open subset parametrizing admissible $G$-covers in which the target curve is smooth.

\begin{proposition}\label{prop:normal_crossings}
    The space 
    $\ov{\cM}^G_{g,S}=\coprod_\rho \ov{\cM}^G_{g,S}(\rho)$ is a normal crossings compactification of ${\cM}^G_{g,S}=\coprod_\rho{\cM}^G_{g,S}(\rho)$.
\end{proposition}
\begin{proof}
   This follows from the fact that $\Adm^\circ_{g,n}(G) \subset \Adm_{g,n}(G)$ is a normal crossings compactification, by the proof of \cite[\S3.23]{mochizuki-geometry}, and $\ov{\cM}^G_{g,S}$ is \'etale over $\Adm_{g,S}(G)$  (Proposition~\ref{prop:etale}).
\end{proof}

\section{Boundary complexes of pointed admissible $G$-covers}\label{section:dual_complex}

In this section we write down the boundary complex for the normal crossings compactification \begin{equation}\label{eq:M0SGrho}{\cM}^G_{0,S}(\rho)\subset \ov{\cM}^G_{0,S}(\rho)\end{equation} when $G$ is prime cyclic (Theorem~\ref{thm:boundary-complex}). This will be used in Section 4, to provide a normal crossings compactification of $\mathcal{H}_{g,n}$ and obtain its boundary complex.
% Just as the boundary complex of $\ov \cM_{g,n}$ is governed by graphs, 
The boundary complex %of $\ov \cM^G_{g,n}$ 
is governed by {\em graph-theoretic admissible covers of graphs}, which we develop below in \S \ref{sec:graph covers}.  

The basic notion of an admissible cover in tropical geometry was established in \cite{caporaso-gonality} and \cite{cavalieri-markwig-ranganathan-admissible}, and hyperelliptic graphs and tropical curves were studied in \cite{bn-hyp} and \cite{hyp}. 
In recent work of Len, Ulircsch, and Zakharov, \cite{LUZ}, the authors classify harmonic $G$-covers of a tropical curve for abelian $G$.
More closely related to this paper, the combinatorics of the stratification of admissible covers spaces by dual graphs is in \cite[\S7]{bertin-romagny-champs}.  Building on this, Schmitt--van Zelm define {\em admissible $G$-graphs}, which are the graphs with $G$-action arising as dual graphs to admissible $G$-covers. They note that a stratum corresponding to an admissible $G$-graph may be disconnected or empty. They also compute the degree of the map from such a stratum to the moduli space of target curves \cite[\S3]{schmitt-van-zelm-intersections}. Implementations in SageMath are available in the package \texttt{admcycles} \cite{delecroix-schmitt-van-zelm-admcycles}.
Closely related, the notion of a {\em graph $G$-cover} associated to a admissible $G$-cover was developed by Galeotti~\cite{galeotti-birational, galeotti-moduli}---see especially~\cite[\S3.1]{galeotti-moduli}---for the purpose of studying the birational geometry, and  singularities, of (coarse spaces of) moduli spaces of genus $g$ curves with a principal $G$-bundle. Our definition is a version of these, undertaken in a case when it becomes possible to explicitly determine the combinatorics of the {\em connected} strata of the boundary.  In other words, by putting into place our restrictions on $g$ and $G$, we are able to give a completely explicit description of the boundary complex of~\eqref{eq:M0SGrho}, which is likely hard in general. See Remarks~\ref{rem:general-case} and~\ref{rem:cmr} for further comments on the general case and for further discussion of the surrounding literature.
%With the abelian hypothesis on $G$, we its group operation is henceforth written additively and the identity element is denoted $0$.

%Given a pointed $G$-admissible cover $(\cP, \wt p_1,\ldots, \wt p_n) \to (\cC, p_1,\ldots,p_n)$ of nodal curves, it is possible to associate some combinatorial data: a morphism of dual graphs $P\to C$, an action $G\times P\to P$ leaving $P\to C$ fixed, and elements $g_1,\ldots,g_n$ of $G$ recording the monodromy around $\wt p_1,\ldots, \wt p_n$ respectively.  But which of these are realized by some pointed $G$-admissible cover? And what are the connected components of the stratification of $\ov{\cM}_{g,n}^G$?  These questions seem delicate, but interesting, in general.  For our purposes, we impose the following conditions:
%\begin{center}
%Take $g=0$, and let $G$ be an abelian group.
%\end{center}
%In this situation both questions can be completely resolved.  

\subsection{Categories of covers of graphs}\label{sec:graph covers} Throughout Section~\ref{section:dual_complex}, let $G$ be a finite prime cyclic group.  In this section, we describe the boundary strata of the compactification \[\cM_{0, S}^{G}(\rho) \hookrightarrow \overline{\cM}_{0, S}^G(\rho),\] showing in Theorem \ref{thm:boundary-complex} that they are in correspondence with \textit{graph-theoretic admissible $G$-covers}, which we will now define. 

A {\em graph} $C = (V,H,i_C,r_C)$ is the data of two finite sets of {\em vertices} $V=V(C)$, and {\em half-edges} $H=H(C)$, together with maps \[ i_C \col H\to H,\qquad r_C \col H\to V\]such that $i_C$ is an involution.  We abbreviate $i= i_C$ and $r = r_C$.  We permit $i$ to have fixed points, and let $L = L(C)$ denote the set of fixed elements of $i$, called {\em legs}.   View $r_C$ as the map taking a half-edge to its incident vertex.
The edge set $E = E(C)$ is the set of pairs $\{h, i(h)\}$ for $i(h)\ne h$; view $i_C$ as the ``other half'' map on the half-edges.    

A morphism of graphs $f\col C\to C'$ is 
given by set maps $f_V \col V\to V'$ and $f_H\col H\to H'$ such that the relevant squares commute:
\[\squarediagramlabel{H}{H}{H'}{H'}{i_C}{f_H}{f_H}{i_{C'}}\qquad \squarediagramlabel{H}{V}{H'}{V'.}{r_C}{f_H}{f_V}{r_{C'}}\]
For a finite set $S$, an {\em $S$-marking} of $C$ is an injection $m = m_C\col S\to L(C)$. It will be convenient {\em not} to require that $m$ is a bijection.  A morphism of $S$-marked graphs $(C,m_C) \to (C', m_{C'})$ is a morphism of graphs $f\col C\to C'$ that preserves the $S$-marking, i.e., $f_H \circ m_C = m_{C'}$.  %An $S$-marked stable graph the dual graph of an $S$-marked stable curve. Precisely: it is a connected $S$-marked graph $C$ and a function $w\col V(G)\to \ZZ_{\ge 0}$ satisfying
%\[2w(v) - 2 + |m^{-1}(v)| > 0\]
%for every vertex $v$.  Its genus is \[|E(C)| - |V(C)| + 1 + \sum_v w(v).\]

\defnow{\label{def:GraphicCovers}Let $G$ be a finite prime cyclic group, $S$ a finite set.
An {\em $S$-marked, admissible $G$-cover of graphs} in genus $0$ is
\enumnow{
\item A morphism $f\col P\to C$ of $S$-marked graphs, such that $C$ is a \textit{stable} $S$-marked tree: for each vertex $v \in V(C)$, we have $|r_C^{-1}(v)| \geq 3$, and that $m_c$ is a bijection between $S$ and the legs of $C$.
\item A left action $\Phi\col G\times P\to P$ leaving $P\to C$ invariant, such that $P\to P/G$ is canonically isomorphic to $P\to C$.
\item A ``monodromy marking'' $\mu\col H(C) \to G$.  Thus every half edge (including legs) of $C$ is assigned an element of $G$.  If $i(h) \neq h$, we require that $\mu(i(h)) = \mu(h)^{-1}$.
\item A function $g\colon  V(P) \to \ZZ_{\geq 0}$; we call $g(v)$ the \textit{weight} or \textit{genus} of $v$.
}
The above data must satisfy:
\begin{enumerate}[(a)]
\item \label{cond:v}For every $v\in V(C)$, $f^{-1}(v) \cong G/\langle \mu(h):  h\in r^{-1}(v)\rangle$ as left $G$-sets, and 
%\[ \sum_{h\in r^{-1}(v)} \mu(h) = 0.\]
% in product notation:
\[ \prod_{h\in r^{-1}(v)} \mu(h) = 1.\]
\item \label{cond:h}For every $h \in H(C)$, $f^{-1}(h) \cong G/\langle \mu(h)\rangle$ as left $G$-sets.
\item{(local Riemann-Hurwitz)}\label{cond:g} For all $v\in V(P)$, writing $w = f(v)$ and $n_w = r_C^{-1}(w)$, the genus $g(v)$ of $v$ is given by
\[2 - 2g(v) =  |\langle \mu(n_w) \rangle|\left(2 - \sum_{h \in n_w} \frac{|\langle\mu(h) \rangle| - 1}{|\langle\mu(h) \rangle|}\right). \]
\end{enumerate}
We will use the boldface notation $\mathbf{P} \to \mathbf{C}$ to indicate a graph-theoretic admissible $G$-cover, with the understanding that this includes all of the data above. When we need to refer to the marking functions, we will write $m_P$ for the marking of $P$ and $m_C$ for the marking of $C$. 
}

It is clear from condition \ref{cond:g} that the genus function $g$ is determined by the monodromy marking $\mu$ as well as the morphism $P \to C$. Moreover, since $C$ is a tree, the data of $C$ and $\mu$, without the $S$-marking, actually determine $P$ and $\Phi$ up to isomorphism. On the other hand, the $S$-marking on $P$ is not in general determined by the $S$-marking on $C$.

\begin{figure}[h]
    \centering
    \includegraphics[scale=1]{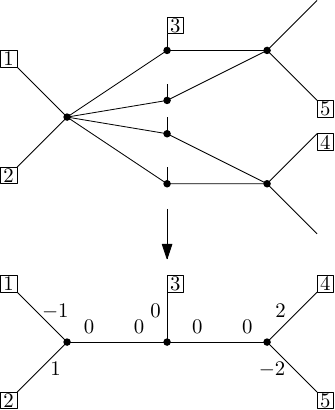}
    \caption{A $G$-cover of $5$-marked graphs, for $G=\ZZ/4\ZZ = \{0,1,2,3\}$. The labels of legs are boxed to avoid confusion with the monodromy marking $\mu: H(C) \to \ZZ/4\ZZ$.}
    \label{fig:Z4cover}
\end{figure}

If $\cP\to \cC$ is an $S$-marked admissible $G$-cover of nodal curves, with $\cC$ a stable $S$-marked curve of genus $0$, then we obtain a corresponding $S$-marked admissible $G$-cover of dual graphs $\mathbf{P}\to \mathbf{C}$.  The meaning of condition \ref{cond:v} is that the subgroup of $G$ stabilizing the generic point of an irreducible component of $\cP$ above a given irreducible component $\cC_v$ of $\cC$ is exactly the subgroup of $G$ generated by the monodromy elements around the special points (nodes and marked points) on $\cC_v$.  Thus, each irreducible component of $\cP$ above $\cC_v$ maps to $\cC_v$ with degree $|\langle \mu(r_C^{-1}(v)) \rangle|$.  The content here is that since $\cC_v$ is rational, $\pi_1(\cC_v)$ is generated by keyhole loops around the special points.  Similarly, the data of a homomorphism $\pi_1(\cC_v)\to G$, for appropriately chosen keyhole loops, is the data of an ordered tuple of elements of $G$ whose product is the identity.  Condition \ref{cond:h} is similar.

\defnow{Let $\mathbf{P} \to \mathbf{C}$ and $\mathbf{P}' \to \mathbf{C}'$ be graph-theoretic $S$-pointed admissible $G$-covers.

\begin{enumerate}
    \item An \textit{isomorphism} $(\mathbf{P} \to \mathbf{C}) \to (\mathbf{P}' \to \mathbf{C}')$ is the data of $G$-equivariant graph isomorphisms $\phi\colon  P \to P'$ and $\psi\colon  C \to C'$, compatible with the marking functions $m_P$ and $m_C$, as well as the monodromy marking $\mu$, which fit into a commutative square.
    \item Let $e \in E(C)$ be an edge. The \textit{edge-contraction} of $\mathbf{P} \to \mathbf{C}$, denoted $(\mathbf{P} \to \mathbf{C})/e$, is obtained by contracting the edge $e$ in $C$, together with its preimages in $P$. The new monodromy marking is obtained by restricting the previous one.
\end{enumerate}
}

\defnow{\label{def:Gamma-category}We write $\Gamma_{0, S}^{G}$ for the category of all graph-theoretic $S$-pointed admissible $G$-covers, where morphisms are given by compositions of isomorphisms and edge-contractions. Given a function $\rho\colon  S \to G$, we put $\Gamma_{0, S}^G(\rho)$ for the full subcategory of $\Gamma_{0, S}^{G}$ on those graph-theoretic $S$-pointed admissible $G$-covers $\mathbf{P} \to \mathbf{C}$ such that the monodromy marking on $\mathbf{C}$ extends $\rho$. Precisely, $\rho = \mu|_{L(C)} \circ m_C$ where $m_C \col S\to L(C)$ is the $S$-marking on $C$.
}

\subsection{The dual complex of the boundary} We now state Theorem~\ref{thm:boundary-complex} on the boundary complex of the space of pointed admissible covers.   
Recall the category of {\em symmetric $\Delta$-complexes} (see \cite{cgp-graph-homology}), i.e., the category $\mathrm{Fun}(\mathsf{FI}^{\mathrm{op}}, \mathsf{Set})$, where $\mathsf{FI}$ is the category of finite sets with injections. 
For $q\ge -1$ an integer, we henceforth write
\[[q] = \{0,\ldots,q\}.\]
This notational convention includes the special case $[-1] = \emptyset$. Given $X\col \mathsf{FI}^{\mathrm{op}}\to\mathsf{Set}$ and an integer $q\ge -1$, write
\[X_q = X([q])\]
for the set of {\em $q$-simplices} of $X$.
%Briefly, a symmetric $\Delta$-complex is a functor $I^{\mathrm{op}}\to\mathsf{Set}$, where $I$ denotes the category of finite sets and injections.  A normal crossings compactification 
\begin{definition}\label{defn:delta_0n(S)}
 Fix $g=0$ and $G$ prime cyclic.  For data $G,S,$ and $\rho$ as above, we define a symmetric $\Delta$-complex
 \[\Delta_{0,S}^G(\rho)\col \mathsf{FI}^\mathrm{op}\to\mathsf{Set}\]as follows.  
 
 For each $q\ge -1$, the set $\Delta_{0,S}^G(\rho)_q$ is the set of isomorphism classes of pairs $(\mathbf{P} \to \mathbf{C}, \omega)$, where
\enumnow{
\item $\mathbf{P} \to \mathbf{C}$ is an object of $\Gamma_{0, S}^{G}(\rho)$
\item  $\omega\col [q]\to E(C)$
is a  bijection, called an {\em edge-labelling}.
}
An isomorphism of pairs $(\mathbf{P} \to \mathbf{C}, \omega) \to (\mathbf{P}' \to \mathbf{C}', \omega')$ is an isomorphism $(\mathbf{P} \to \mathbf{C}) \to (\mathbf{P}' \to \mathbf{C}')$ such that if $\psi: \mathbf{C} \to \mathbf{C}'$ is the induced isomorphism on targets, we have $\omega' = \psi \circ \omega$ as maps $[q] \to E(\mathbf{C}')$.

For morphisms, given $i\col [q']\hookrightarrow [q]$, and given a graph-theoretic admissible cover $\mathbf{P}\to \mathbf{C}$ as above, contract the edges $E(C) - \omega(i([q']))$, to obtain a new object of $\Gamma_{0, S}^{G}(\rho)$, and take the unique edge-labelling by $[q']$ which preserves the order of the remaining edges.
\end{definition}

\thmnow{\label{thm:boundary-complex} Let $G$ be a prime cyclic group, $S$ a finite set.  There is an isomorphism of symmetric $\Delta$-complexes
\[\Delta_{0,S}^G(\rho)\cong \Delta({\cM}^G_{0,S}(\rho)\subset \ov{\cM}^G_{0,S}(\rho)).\]
}
\begin{proof}
Let us start with the stratification of the boundary of $\Adm_{0,S}(G;\rho)$. The space $\Adm_{0,S}(G;\rho)$ is nonempty if and only if $\prod_{s\in S} \rho(s) = 1_G$. The boundary complex of $\Adm^\circ_{0,S}(G;\rho)\subset\Adm_{0,S}(G;\rho)$ is the complex of trees $C$ with a bijective $S$-marking $m\col S\to L(C)$, together with a monodromy marking $\mu\col H(C)\to G$ extending $\rho$, which must satisfy, for every vertex $v\in V(C)$ and $e = \{h_1,h_2\}\in E(C)$, \[\prod_{h\in r^{-1}(v)} \mu(h) = 1, \qquad \mu(h_1)\mu(h_2)=1.\]  
The stratum of the boundary indexed by such a triple $(C,m,\mu)$ is indeed connected, since it is, up to finite quotient, isomorphic to a product $\prod_{v\in V(C)} \Adm_{0,n_v}(G;\mu_v)$ of varieties that are themselves connected, see Equation \eqref{eq:adm-gi}.
More formally, as a symmetric $\Delta$-complex, the boundary complex has a $q$-simplex for every such datum $(C, m, \mu)$ together with an arbitrary bijective edge-labelling $\omega\col [q] \to E(C)$, one for each isomorphism class of $(C,m,\mu, \omega)$.  

%Now since $\ov{\cM}_{0,S}^G(\rho)$ is \'etale over $\Adm_{0,S}(G;\rho)$, there is a morphism of boundary complexes from that of the former to that of the latter.  We now study the fibers of this morphism.  
Suppose 
\[(C, m\col S\to L(C), \mu\col H(C)\to G)\]
is a stable $S$-marked tree with monodromy marking $\mu$ as above. 
For $v\in C$,  write $n_v = r^{-1}(v)$  for the set of half-edges (including legs) at $v$, and write \[G_v  = \langle \mu(h)\col h\in n_v\rangle.\]  
Let $\mu_v$ be the restriction of $\mu$ to $n_v$.    As noted above, $(C,m,\mu)$ indexes a boundary stratum of  $\Adm_{0,S}(G;\rho)$.  
The preimage in $\ov{\cM}_{0,S}^G(\rho)$ of this stratum is isomorphic to the variety
\begin{equation}\label{eq:strata-upstairs}  \!\! \prod_{v\in V(C)} \!\! \left( \cM_{0,n_v}^G(\mu_v) \,/\,G^{E(C)} \right )\end{equation}
e.g., by \cite[\S2]{petersen-operad}.  Let us explain the action of $G^{E(C)}$ in~\eqref{eq:strata-upstairs}. For a given edge $e = \{h,h'\}$, incident to vertices $v$ and $v'$, the copy of $G$ indexed by $e$ acts by translating the lifted marked point indexed by $h$, respectively $h'$, in the moduli space $\cM_{0,n_v}^G(\mu_v)$, respectively $\cM_{0,n_{v'}}^G(\mu_{v'}).$  (In general, $G$ would also change the values of the marking functions $\mu_v(h)$ and $\mu_{v'}(h')$, respectively, by conjugation, but $G$ is abelian here.)

The variety~\eqref{eq:strata-upstairs} may not be connected, and it remains to  describe
%The only loose end in our description of the boundary of $\ov{\cM}_{0,S}(G;\rho)$ is a description of 
its  connected components. %since in general they are not connected.  
For each $v\in V(C)$, let \[X_v = \{\mathrm{Fun}(n_v, G/G_v)\}/G\] where the quotient is with respect to the $G$-action on $G/G_v$.  From Proposition~\ref{prop:cc-smooth-case}, the connected components of ~\eqref{eq:strata-upstairs} are in bijection with
\begin{equation}\label{eq:cut-up-graph} \left(\prod_{v\in V(C)} X_v \right)\!/ G^{E(C)}.\end{equation}

The last step is a combinatorial identification of~\eqref{eq:cut-up-graph} with the set of isomorphism classes of graph-theoretic $S$-pointed admissible $G$-covers.  Let us begin by considering local data at a single vertex $v\in V(C)$. Consider an element $f_v\in X_v$, together with the data of $\mu|_{n_v} \col n_v \to G$.  From $f_v$ and $\mu|_{n_v}$ we can extract a graph-theoretic $n_v$-pointed admissible cover involving graphs with legs but no edges: $C_v$ is a single vertex, with legs $n_v$; $V(P_v) = G/G_v$ as a left $G$-set; and above each leg $h\in n_v$ of $C$ is a set of legs in $P_v$ isomorphic to $G/\langle \mu(h)\rangle,$ with root map compatible with the map $G/\langle \mu(h)\rangle \to G/G_v$.  Finally, $P_v$ has $S$-marking given by $f_v$.  

Continue to fix a stable $S$-marked tree $C$ and monodromy marking $\mu$ on $C$.  Now, given $(f_v)_v \in \prod X_v$, we assemble the local picture above into an admissible cover of graphs. For every edge $e = \{h,h'\}$ of $C$, with root vertices $v=r(h)$ and $v'=r(h')$, the half-edges of $P_v$ above $h$ and the half-edges of $P_{v'} $ above $h'$ are each isomorphic to $G/\langle \mu(h) \rangle = G/\langle \mu(h') \rangle $ as $G$-sets. There is a unique $G$-equivariant bijection between these two sets that sends the chosen lift of $h$ to the chosen lift of $h'$, and another choice of lifts of $h$ and $h'$ produce the same bijection if they are related to the original choices by the same element of $G$.  Therefore these identifications glue the half-edges above $h$ and $h'$ into edges above $e$, obtaining a graph-theoretic admissible cover $P\to C$ which was independent of the action of $G^{E(C)}$.  It is straightforward to reverse this process, giving an element of the set~\eqref{eq:cut-up-graph} starting from a graph-theoretic admissible cover.  
\end{proof}

\begin{remark}\label{rem:general-case}
Theorem~\ref{thm:boundary-complex} furnishes an explicit description of the symmetric $\Delta$-complex 
\begin{equation}\label{eq:the-boundary-complex}\Delta(\cM_{g,n}^G \subset \overline{\cM}_{g,n}^G) \end{equation}
when $g=0$ and $G$ is prime cyclic.  It is sufficiently explicit that it can be programmed, and indeed we carry out computer calculations for the results in Appendix~\ref{sec:calculations}.  Without restrictions on $G$ and $g$, it is still possible to give a general description of~\eqref{eq:the-boundary-complex} using the framework of {\em graphs of groups}, roughly, decorating vertices of graphs with fundamental groups of punctured curves.  This idea will appear in future work by M.~Talpo, M.~Ulirsch, and D.~Zakharov; we thank Ulirsch for bringing it to our attention.  This general description is not explicit in the above sense. It involves the very interesting sub-question of determining the connected components of the spaces $\cM_{g,n}^G$ in general; compare with Proposition~\ref{prop:adm-m0n-trivial}.  We also refer to forthcoming work of P.~Souza, that constructs~\eqref{eq:the-boundary-complex} in the case of $G$ cyclic with $g$ arbitrary, and identifies it as the nonarchimedean skeleton of the toroidal pair.  Moreover, that work is a precursor to further work by Y.~El Maazouz, P.~Helminck, F.~R\"ohrle, P.~Souza, and C.~Yun studying the homotopy type of boundary complexes of unramified $\ZZ/p\ZZ$ covers for $g=2$.
\end{remark}

\begin{remark}\label{rem:cmr}
The graph-theoretic admissible $G$-covers in this paper (Definition~\ref{def:GraphicCovers}) are exactly what are needed for a precise description of the boundary complex (Theorem~\ref{thm:boundary-complex}).  Thus they are reasonably expected to be similar to, but distinct from, the spaces of covers of tropical curves appearing in \cite{caporaso-gonality}, in \cite{cavalieri-markwig-ranganathan-admissible}, and the references therein. 
The work %of Cavalieri--Markwig--Ranganathan's on tropicalizations of the space of admissible covers 
\cite{cavalieri-markwig-ranganathan-admissible} on tropicalizations of the space of admissible covers is 
an important comparison point for this paper.  Rather than $G$-covers, they  study the admissible covers compactification of the Hurwitz space of degree $d$ covers of smooth curves with fixed target genus $h$ and fixed ramification profiles (and hence fixed source genus $g$) over $n$ marked branch points in the target.  All of the inverse images of the branch points are also marked.  This moduli space is %the normalization of the Harris-Mumford admissible covers compactification, and is 
canonically isomorphic to a cover of a component of the space $\Adm_{h,n}(S_d)$. % of Abramovich-Corti-Vistoli admissible covers. 
%(Specifically, the component is specified by the ramification data, and the cover corresponds to the extra data of the source markings.)  
In \cite{cavalieri-markwig-ranganathan-admissible} the boundary complex, which may be identified with the link of the skeleton of the Berkovich analytification \cite{acp}, is {\em compared}, but not identified, with a certain space of tropical admissible covers, via a surjective morphism of generalized cone complexes from the former to the latter. The failure of this surjection to be an isomorphism is due to multiplicities fully accounted for in~\cite[\S4.2.4]{cavalieri-markwig-ranganathan-admissible}, and is related to Remark~\ref{rem:general-case} above.

\end{remark}
%Let $X$ be the boundary complex ${\cM}_{0,2g+2+n}^{\ZZ/2\ZZ}(0^n, 1^{2g+2}) \subset \ov{\cM}_{0,2g+2+n}^{\ZZ/2\ZZ}(0^n, 1^{2g+2})$, regarded as a symmetric $\Delta$-complex \cite{cgp-graph-homology}.
%\propnow{The set of $q$-simplices are isomorphism classes of data
%\[(C\xra{\alpha} T, W, m, \tau)\]
%where
%\itemnow{\item $C\to T$ is a harmonic double cover of graphs in which $T$ is a tree; let $j \col C\to C$ denote the hyperelliptic involution on $C$;
%\item $W\col \{1,\ldots,2g+2\}\to V(C)^j$ and  $m\col \{1,\ldots,n\} \to V(C)$ are functions,
%\item $\tau\col \{0,\ldots,q\}\xra{\cong} E(T)$ is any bijection
%}
%such that 
%\itemnow{\item the induced maps $\alpha \circ W$ and %$\alpha \circ m$ together make $T$ a stable tree.}
%}

%%%%%%

\section{Compactifications of $\cH_{g,n}$}
\label{section:compactifications}

Let $g\ge 2$ and $n\ge 0$. Throughout this section we will fix
\[S = \{1, \ldots, n\} \cup \{w_1, \ldots, w_{2g + 2} \}, \]
and fix $G = \ZZ/2\ZZ = \{0,1\}$.  We also define $\rho\colon  S \to \ZZ/2\ZZ $ by $\rho(i) = 0$ for all $i \in \{1, \ldots, n\}$, and $\rho(w_k) = 1$ for $k \in \{1, \ldots, 2g+2\}$. We will discuss how the stack quotient
\[[\ov{\cM}^{\ZZ/2\ZZ}_{0,S}(\rho)/S_{2g+2}]\]provides a normal crossings compactification of $\cH_{g, n}$, and give an explicit description of the dual complex $\Theta_{g, n}$ of this compactification. The description will be in terms of the dual complexes studied in the previous section. We first consider the case of labelled Weierstrass points, and then quotient out by $S_{2g+2}$.

\subsection{The complex $\wt{\Theta}_{g, n}$}\label{section:tildeH_gn}
First, let $\wt{\cH}_{g,n}$ denote the moduli stack of hyperelliptic curves of genus $g$ with $n$ distinct marked points %$p_1,\ldots, p_{n}$, 
and $2g+2$ labelled Weierstrass points.  The symmetric group on $2g+2$ letters permutes the labels on Weierstrass points, and
\[\cH_{g,n} \cong [\wt{\cH}_{g,n}/S_{2g+2}].\]
In this subsection, we will provide a normal crossings compactification of $\wt{\cH}_{g,n}$ and give the corresponding dual complex. Then, we will quotient out by $S_{2g+2}$ to give a normal crossings compactification of $\cH_{g,n}$. 

\begin{figure}[h]
    \centering
    \includegraphics[scale=0.9]{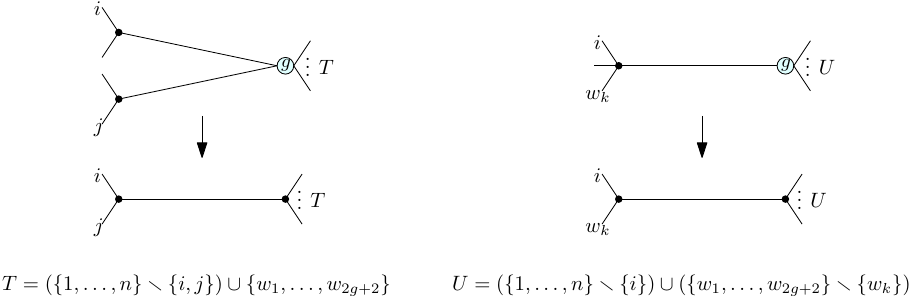}
    \caption{Two graph-theoretic admissible $G$-covers, where $G = \ZZ/2\ZZ = \{0,1\}.$}
    \label{fig:two-covers}
\end{figure}

 In $\wt{\cH}_{g,n}$, a marked point {\em is} allowed to coincide with a Weierstrass point, and two marked points are allowed to form a conjugate pair under the hyperelliptic involution. Because of this,
two types of graphs will require special attention.
\begin{definition}\label{defn:orderedH_gn} 
We call the following graph-theoretic admissible covers type (1) and type (2) respectively:
\enumnow{
\item For distinct $i,j\in\{1,\ldots,n\}$, the admissible cover of graphs in Figure~\ref{fig:two-covers} on the left.  
\item For each $i\in \{1,\ldots,n\}$ and $w_k \in \{w_1, \ldots, w_{2g + 2}\}$, the admissible cover of graphs in Figure~\ref{fig:two-covers} on the right.
}
\end{definition}

\propnow{There is an open inclusion \[\wt{\cH}_{g,n} \hookrightarrow \ov{\cM}^{\ZZ/2\ZZ}_{0,S}(\rho)\] which is a normal crossings compactification, and whose boundary complex $\wt{\Theta}_{g, n}$ is isomorphic to the subcomplex of 
\[\Delta_{0,S}^{\ZZ/2\ZZ}(\rho)\]
on simplices whose vertices are not of type (1) or (2) in Definition \ref{defn:orderedH_gn}. 
\label{cor:htilde}
}

\begin{proof}
Let $\wt{\cH}^\circ_{g,n}$ denote the open substack of  $\widetilde{\cH}_{g,n}$ in which a marked point may not collide with a Weierstrass point, and two marked points may not form a conjugate pair. Then
\[\wt{\cH}^\circ_{g,n} \cong \cM_{0,S}^{\ZZ/2\ZZ}(\rho),\]
where $\cM_{0,S}^{\ZZ/2\ZZ}(\rho)$ denotes the interior of the moduli space $\ov{\cM}_{0,S}^{\ZZ/2\ZZ}(\rho)$ of pointed admissible covers. 
 We define a partial compactification $\cH_{g, n}^*$ of $\wt{\cH}_{g, n}^{\circ}$, such that
 \[\wt{\cH}_{g, n}^{\circ} \subset \cH_{g, n}^* \subset \ov{\cM}^{\ZZ/2\ZZ}_{0,S}(\rho), \]
 and the second inclusion is normal crossings.
In $\ov{\cM}^{\ZZ/2\ZZ}_{0,S}(\rho)$, define $\cH_{g,n}^*$ to be the open complement of all boundary divisors except for those corresponding to dual graphs of type (1) or (2) (see Definition \ref{defn:orderedH_gn}). Since $\cH_{g, n}^*$ is the complement of a subset of the boundary divisors, the divisor
\[\ov{\cM}_{0, S}^{\ZZ/2\ZZ}(\rho) \smallsetminus \cH_{g, n}^* \] still has normal crossings. Stabilization gives a canonical isomorphism $\cH_{g,n}^* \cong \wt{\cH}_{g,n}$ which is equivariant with respect to the action of $S_n$, thus giving the first part of the result.

We now turn our attention to the boundary complex. 
Denote by $\Delta_{0,S}^{\ZZ/2\ZZ}(\rho)$ the dual complex of the compactification \[\wt{\cH}^\circ_{g,n} \cong \cM_{0,S}^{\ZZ/2\ZZ}(\rho)\subset \ov{\cM}_{0,S}^{\ZZ/2\ZZ}(\rho).\]
The target graphs of type (1) and (2) in Definition \ref{defn:orderedH_gn} have one edge, and correspond to vertices in $\Delta_{0,S}^{\ZZ/2\ZZ}(\rho)$.
Then, the boundary complex  $\wt{\Theta}_{g, n}$ of the inclusion
\[\wt{\cH}_{g,n} \subset \ov{\cM}^{\ZZ/2\ZZ}_{0,S}(\rho)\] is the subcomplex of 
$\Delta_{0,S}^{\ZZ/2\ZZ}(\rho)$
determined by those simplices which have no vertices of type (1) or (2) in Definition~\ref{defn:orderedH_gn}.   
\end{proof}

 Let us now describe the complex $\Delta_{0,S}^{\ZZ/2\ZZ}(\rho)$ in more detail.  
 %Write $\{0,1\}$ for the elements of $\ZZ/2\ZZ$.
 Its $q$-simplices  are given by isomorphism classes of pairs $(\mathbf{P} \to \mathbf{C}, \omega)$, where $\mathbf{P} \to \mathbf{C}$ is an object of the category $\Gamma_{0, S}^{\ZZ/2\ZZ}(\rho)$ (Definition~\ref{def:Gamma-category}), and $\omega\colon  [q]\to E(C)$ is an edge-labelling. Moreover, on $L(C)$, the monodromy marking $\mu$ satisfies $\mu(m_C(j)) = 0$ if $j \in \{1, \ldots, n\}$, and $\mu(m_C(j)) = 1$ if $j \in \{w_1, \ldots, w_{2g + 2}\}$. We will call the elements of \[m_C(\{w_1, \ldots, w_{2g + 2} \}) \subset L(C) \]
the \textit{branch legs} of $C$.

Notice that the above conditions on $\mu|_{L(C)}$ suffice to determine $\mu$ on all other half-edges of $C$, by condition (1) of Definition \ref{defn:delta_0n(S)}. Call a vertex $v \in V(C)$ a \textit{leaf} vertex if it is incident to only one edge. If a leaf vertex $v \in V(C)$ supports an odd number of branch legs, then the non-leg half edge $h$ incident to $v$ must satisfy $\mu(h) = 1$. On the other hand, if a leaf vertex $v$ supports an even number of branch legs, then the non-leg half edge $h$ incident to $v$ must satisfy $\mu(h) = 0$. Proceeding inductively, this determines $\mu$ on all half-edges incident to non-leaf vertices of $C$ as well.

\begin{figure}[h]
    \centering
\includegraphics[scale=1.25]{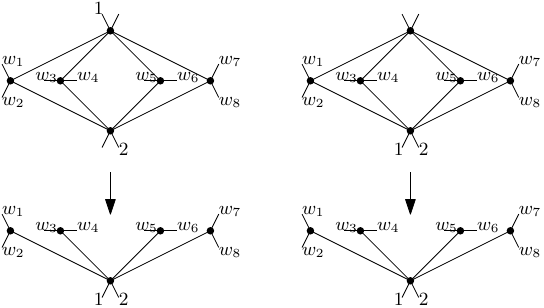}
    \caption{A $\{1, 2\} \cup \{w_1, \ldots, w_8\}$-marked stable tree $C$, together with the two lifts of $m_C$ to a marking $m_P$. These non-isomorphic lifts are determined by a choice of element in the fiber over each leg marked by $\{1, 2\}$ on $C$, and two such choices define the same graph-theoretic admissible $\ZZ/2\ZZ$-cover if they differ by the $\ZZ/2\ZZ$-action on $P$.}
    \label{marking_example}
\end{figure}

This discussion implies that given the monodromy data $\rho$ and an $S$-marked stable tree $C$, the only additional data required to determine an object of the category $\Gamma_{0, S}^{\ZZ/2\ZZ}(\rho)$ is 
a lift of the marking function 
%\[m_C : \{1, \ldots, n \} \cup \{w_1 \ldots, w_{2g + 2}\} \to L(C) \] to a function
%\[m_P: \{1, \ldots, n \} \cup \{w_1 \ldots, w_{2g + 2}\} \to L(P) \]
$m_C\col S\to L(C)$ to a function $m_P\col S\to L(P)$ 
such that the diagram
%\[ \begin{tikzcd}
%& & L(P) \arrow[dd]\\
%& \{1, \ldots, n \} \cup \{w_1 \ldots, %w_{2g + 2}\} \arrow[ur, "m_P"] \arrow[dr, "m_C"] &\\
%& & L(C) 
%\end{tikzcd} \]
\[{\xymatrix@R=6mm@C=12mm{
& L(P)\ar[d] \\ S \ar[ur]^{m_P} \ar[r]_{m_C} &  L(C)
}}\]
commutes.
(Note that the morphism of graphs $P \to C$, without the marking function on $P$, is already determined by $C$ and $\mu$.) Moreover, since each branch leg in $C$ has a unique preimage in $P$, one only needs to choose, for each $i \in \{1, \ldots, n\}$, a leg in the preimage of $m(i) \in L(C)$. Two such choices are equivalent if they differ by the $\ZZ/2\ZZ$-action on $P$. See Figure \ref{marking_example} for an example.

 % with prescribed monodromy at marked points \cite{abramovich-corti-vistoli-twisted, jarvis-kaufmann-kimura-pointed}.  

 \subsection{The complex $\Theta_{g, n}$}
We now construct a normal crossings compactification of $\cH_{g,n}$ and the corresponding dual complex $\Theta_{g, n}$.

By Proposition \ref{cor:htilde}, in order to pass from $\Delta_{0, S}^{\ZZ/2\ZZ}(\rho)$ to $\wt{\Theta}_{g, n}$, we remove all edge-labelled pairs $(\mathbf{P} \to \mathbf{C}, \omega)$ such that $\mathbf{P} \to \mathbf{C}$ admits a contraction to covers of type (1) or (2) in Definition \ref{defn:orderedH_gn}. To that end, let \[\Gamma_{0, S}^{\ZZ/2\ZZ, *}(\rho)\]  be the full subcategory of $\Gamma_{0, S}^{\ZZ/2\ZZ}(\rho)$ on those covers which do not admit a contraction to covers of type (1) or (2). 
\begin{definition}
We define the category $\Gamma^{\cH}_{g, n}$ as follows.
\begin{enumerate}
    \item The objects are $S_{2g + 2}$-orbits of objects of $\Gamma_{0, S}^{\ZZ/2\ZZ, *}(\rho)$. Precisely, the objects are covers $\mathbf{P} \to \mathbf{C}$, where
    \begin{enumerate}[(a)]
    \item $\mathbf{C} = (C, m_C)$ is the data of a stable tree $C$ with $2g + 2 + n$ legs, together with an injective function $m_C\col  \{1, \ldots, n\} \to L(C)$. % such that there are no vertices $v \in V(C)$ such that \[|r^{-1}(v)| = 3,\quad |L(C) \cap r^{-1}(v)| = 2,\quad\text{and}\quad |m_C^{-1}(v)| = 1.\]
    \item $\mathbf{P} = (P, m_P)$, where $P$ is the unique graph-theoretic admissible $\ZZ/2\ZZ$-cover of $C$ obtained by declaring each unmarked leg to have monodromy $1 \in \ZZ/2\ZZ$ and each marked leg to have monodromy $0$, and $m_P \col \{1, \ldots, n\} \to L(P)$ is a marking of $L(P)$ such that $m_P(i)$ is a leg in the inverse image of $m_C(i)$ for all $i$. 
        \end{enumerate}
    The cover $\mathbf{P}\to\mathbf{C}$ is required to satisfy: 
    \begin{itemize}
        \item No contraction to type (1): If $v\in V(C)$ has $|r^{-1}(v)|=3$ and $m_C^{-1}(r^{-1}(v)) = \{i,j\}$, then a single vertex of $P$ supports markings $i$ and $j$---in other words, \[r_P(m_P(i)) = r_P(m_P(j)).\]      
        \item No contraction to type (2): No vertex $v \in V(C)$ satisfies \[|r^{-1}(v)| = 3,\quad |L(C) \cap r^{-1}(v)| = 2,\quad\text{and}\quad |m_C^{-1}(r^{-1}(v))| = 1.\]
    \end{itemize}

    \item The morphisms are compositions of isomorphisms and edge-contractions.
\end{enumerate}
\end{definition}

% this seemed a little repetitive
%We now show that $\cH_{g,n}\subset [\ov{\cM}^{\ZZ/2\ZZ}_{0,S}(\rho)/S_{2g+2}]$ is a normal crossings compactification and 
%provide the boundary complex $\Theta_{g,n}$.

\begin{proposition}\label{prop:Theta_gn}
The inclusion $\cH_{g,n}\subset [\ov{\cM}^{\ZZ/2\ZZ}_{0,S}(\rho)/S_{2g+2}]$ is a normal crossings compactification, and the boundary complex $\Theta_{g,n}$ 
has the following explicit description.
\begin{enumerate}
\item The set of $q$-simplices $\left(\Theta_{g, n}\right)_q$ is the set of isomorphism classes of pairs $(\mathbf{P} \to \mathbf{C}, \omega)$ where $\mathbf{P} \to \mathbf{C}$ is an object of $\Gamma_{g, n}^{\cH}$, and $\omega\col [q] \to E(C)$ is an edge-labelling.
\item Given an injection $\iota\col  [q'] \hookrightarrow [q]$, we define $\iota^*(\mathbf{P} \to \mathbf{C}, \omega) \in \left(\Theta_{g, n}\right)_{q'}$ by contracting those edges which are not in the image of $\iota$, and taking the unique induced edge-labelling which preserves the order of the remaining edges.
\end{enumerate}
\end{proposition}

\begin{proof}
Since the action of $S_{2g+2}$ on $\wt{\cH}_{g,n} \subset \ov{\cM}^{\ZZ/2\ZZ}_{0,S}(\rho)$ preserves $\wt{\cH}_{g,n}$ and sends strata isomorphically to strata, we have that 
$$\cH_{g,n} \cong [\wt{\cH}_{g,n} / S_{2g+2}] \subset [\ov{\cM}^{\ZZ/2\ZZ}_{0,S}(\rho) / S_{2g+2}]$$
is a normal crossings compactification with boundary complex equal to \[\Delta(\wt{\cH}_{g,n} \subset \ov{\cM}^{\ZZ/2\ZZ}_{0,S}(\rho)) / S_{2g+2} = \widetilde{\Theta}_{g, n}/S_{2g + 2},\] and the described symmetric $\Delta$-complex is precisely the quotient of $\widetilde{\Theta}_{g, n}$ by $S_{2g + 2}$.
\end{proof}
% There is a functor
% \[ \Gamma_{0, S}^{\ZZ/2\ZZ, *}(\rho) \to \Gamma^{\cH}_{g, n}, \] given by forgetting the ordering of the branch legs. On objects, this functor is the quotient by the action of $S_{2g + 2}$.
As a direct result of Proposition \ref{prop:Theta_gn}, we have the following corollary identifying the weight zero compactly supported cohomology of $\cH_{g, n}$ with the reduced cohomology of $\Theta_{g, n}$: see \cite[Theorem 5.8]{cgp-graph-homology}.
\begin{corollary}\label{cor:WeightZeroIdentification}
    For each $i$, there are canonical $S_n$-equivariant isomorphisms
    \[W_0 H^i_c(\cH_{g, n};\QQ) \cong \widetilde{H}^{i - 1}(\Theta_{g, n};\QQ) \cong \widetilde{H}_{i - 1}(\Theta_{g, n};\QQ)^{\vee}, \]
    where $\widetilde{H}^*$ and $\widetilde{H}_*$ denote reduced cohomology and homology, respectively.
\end{corollary}

We now establish some conventions for working with objects of the category $\Gamma^{\cH}_{g, n}$. \begin{definition}\label{def:weight} Given an object $\mathbf{P} \to \mathbf{C}$ of $\Gamma^{\cH}_{g, n}$, we define the \textit{weight} of a vertex $v \in V(\mathbf{C})$ to be the number of unmarked legs based at $v$.\end{definition}
\noindent The total weight of the vertices of $C$ is $2g+2$. The weight in this sense should not be confused with the notion of vertex weights corresponding to genera of irreducible curves. The two notions of vertex weight are related by the Riemann-Hurwitz formula.

\begin{figure}
    \centering
    \includegraphics[scale=1.25]{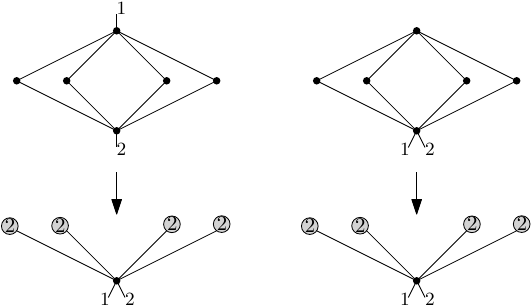}
    \caption{The images of the graph-theoretic admissible $\ZZ/2\ZZ$-covers in $\Gamma_{0, S}^{\ZZ/2\ZZ, *}(\rho)$ from Figure \ref{marking_example}, under the functor $\Gamma_{0, S}^{\ZZ/2\ZZ, *}(\rho) \to \Gamma_{g, n}^{\cH}$. The number of unmarked legs at a vertex of a target tree is indicated by the weight function. We do not depict any unmarked legs of the source graph, since they are determined by the legs of the target.}
    \label{markingquotientfig}
\end{figure}

\begin{figure}
    \centering
    \includegraphics[scale=1.25]{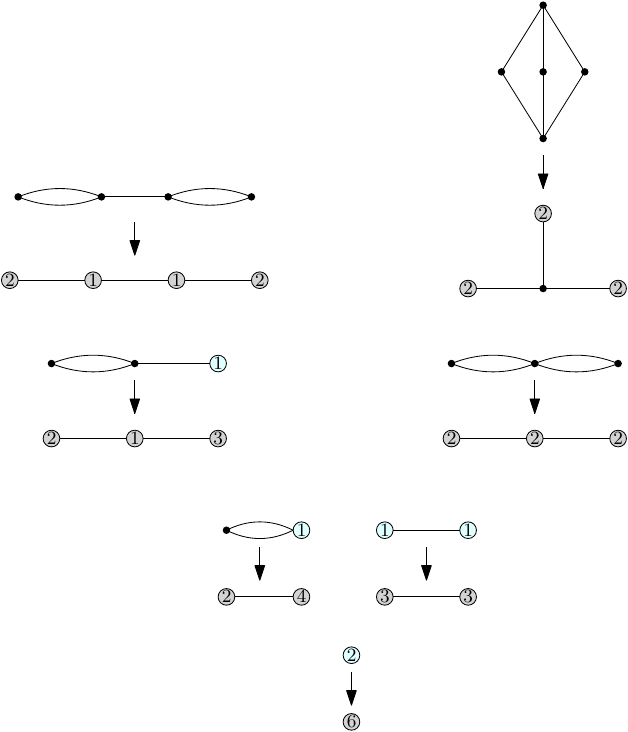}
    \caption{The set of isomorphism classes of $\Gamma^{\cH}_{g, n}$-objects for $g = 2$ and $n = 0$.}
    \label{genus2fig}
\end{figure}

When depicting objects of $\Gamma_{g, n}^{\cH}$, we adopt the following conventions. Instead of drawing the unmarked legs of $\mathbf{C}$, we will label each vertex of $\mathbf{C}$ with its weight. To avoid confusion with the genera of vertices in the source graph, we will depict the weight of a vertex in $\C$ with the color grey, and genera of vertices with blue. Since each unmarked leg of $C$ has a unique preimage in $P$, we will not draw those legs of $P$. When a leg of $C$ has two preimages in $P$, so only one is marked, we will suppress the other leg. See
Figure \ref{markingquotientfig} for the images of the $\Gamma_{0, S}^{\ZZ/2\ZZ, *}(\rho)$ objects from Figure \ref{marking_example} under the functor to $\Gamma^{\cH}_{g, n}$. See Figure \ref{genus2fig} for a complete list of isomorphism classes of $\Gamma_{g, n}^{\cH}$-objects when $g = 2$ and $n = 0$. 

\begin{remark} We remark on the case $n=0$.  In this case, the symmetric $\Delta$-complex $\Theta_{g,0}$ is isomorphic to the quotient of the dual complex
\begin{equation}\label{eqn: dual-complex-m0n}
\Delta_{0, 2g+2}:= \Delta\left(\cM_{0, 2g+ 2} \subset\overline{\cM}_{0, 2g + 2}\right)\end{equation}
by the $S_{2g + 2}$-action permuting the marked points. The dual complex (\ref{eqn: dual-complex-m0n}) is the moduli space of $(2g+2)$-marked tropical curves of genus zero and volume one \cite{cgp-marked}, also known as the space of phylogenetic trees \cite{ardila-klivans-bergman,bhv,robinson-whitehouse-tree}. The identification
\[\Theta_{g,0} = \Delta_{0,2g + 2} /S_{2g + 2} \]
%follows either from the observation that when $n = 0$, the morphism
%\[[\overline{\cM}_{0, S}^{\ZZ/2\ZZ}(\rho)/S_{2g + 2}] \to
% [\overline{\cM}_{0, 2g + 2}/S_{2g + 2}] \]
% sending a cover to its marked target curve is a $\ZZ/2\ZZ$-gerbe, or directly from our description of the category $\Gamma_{g}^\cH$.
can be seen directly from our description of the category $\Gamma_{g}^\cH$, and holds despite the fact that the morphism
\[[\overline{\cM}_{0, 2g+2}^{\ZZ/2\ZZ}(\rho)/S_{2g + 2}] \to
 [\overline{\cM}_{0, 2g + 2}/S_{2g + 2}] \]
is not an isomorphism or even a $\ZZ/2\ZZ$-gerbe, due to the possible presence of extra automorphisms, more than $\ZZ/2\ZZ$, in the source curves of $\ZZ/2\ZZ$-admissible covers.
\end{remark}

\section{Acyclic subcomplexes of $\Theta_{g, n}$}\label{section:acyclic}

In this section we will study the cellular chain complex of $\Theta_{g, n}$, establishing Theorem \ref{AcyclicSubcomplexes} below, which states that several natural subcomplexes are acyclic. This will allow us to prove Proposition \ref{mainthm:smalln} later in this section. The acyclicity results will be used in Section \ref{section: graphsumformula} to obtain Theorem \ref{mainthm:FrobChar}. 
\begin{theorem}\label{AcyclicSubcomplexes} Fix $g \geq 2$ and $n \geq 0$. Then the following subcomplexes of $\Theta_{g, n}$ have vanishing reduced rational homology:
\begin{enumerate}
    \item the \textit{repeated marking locus} $\Theta_{g, n}^{\mathrm{rep}}$, namely the subcomplex determined by those $\Gamma_{g,n}^{\cH}$-objects $\mathbf{P} \to \mathbf{C}$ such that there exists $v \in V(\mathbf{P})$ supporting at least two markings from $\{1, \ldots, n\}$;
    \item the \textit{weight $3$ locus} $\Theta_{g, n}^{\geq 3}$, determined by those $\Gamma_{g,n}^{\cH}$-objects $\mathbf{P} \to \mathbf{C}$ such that $\mathbf{C}$ has a vertex of weight at least $3$ (Definition~\ref{def:weight}); and
    \item the intersection $\Theta_{g, n}^{\mathrm{rep}} \cap \Theta_{g, n}^{\geq 3}$.
\end{enumerate}
\end{theorem}

\remnow{%While we do not need them for this paper, 
There are stronger statements that are also true, namely that the three subspaces of the space $\Theta_{g,n}$ corresponding to (1), (2), and (3) are in fact contractible. It is possible to convert the proofs below, of vanishing reduced rational homology, to proofs of contractibility, using the {\em vertex property} technique of \cite[\S4]{cgp-marked}.}

% Additionally, when $n \leq 1$, the reduced rational homology of $\Theta_{g, n}$ vanishes: this corresponds to the first part of Proposition \ref{mainthm:smalln}, under the identification of the weight zero compactly supported cohomology of $\cH_{g, n}$ with the reduced rational cohomology of $\Theta_{g, n}$.
% \begin{theorem}\label{AcyclicSmallN}
% Fix $g\geq 2$ and $n \leq 1$. Then for all $i$,
% \[ \widetilde{H}_i(\Theta_{g, n};\QQ) = 0. \]
% \end{theorem}

 \subsection{The cellular chain complex of $\Theta_{g, n}$} Following \cite[\S 3]{cgp-graph-homology}, the reduced rational homology of $\Theta_{g, n}$ is computed by the graph complex $\mathcal{C}_*^{(g, n)}$ described as follows. In degree $p$, $\mathcal{C}_p^{(g, n)}$ is spanned by pairs $(\mathbf{P} \to \mathbf{C}, \omega)$ where $\mathbf{P} \to \mathbf{C}$ is an object of $\Gamma^{\mathcal{H}}_{g, n}$, and $\omega\colon   [p] \to E(\mathbf{C})$ is a bijective edge-labelling. These pairs are subject to the relation
$(\mathbf{P} \to \mathbf{C}, \omega) = \mathrm{sgn}(\rho) (\mathbf{P} \to \mathbf{C}, \omega \circ \rho)$ whenever $\rho \in S_{p + 1} = \Aut([p])$.

The differential $\partial\colon  \mathcal{C}_p^{(g, n)} \to \mathcal{C}_{p - 1}^{(g, n)}$ is given by the signed sum of edge contractions:
\[ \partial (\mathbf{P} \to \mathbf{C}, \omega) = \sum_{i \in [p]} (-1)^i (\delta^{i})^{*}(\mathbf{P} \to \mathbf{C}, \omega), \]
where $\delta^i \colon  [p - 1] \to [p]$ is the unique order-preserving injection which misses $i$.

To prove Theorem \ref{AcyclicSubcomplexes}, we will show that the corresponding sub-chain complexes of $\mathcal{C}^{(g, n)}_*$ are acyclic. %To prove Theorem \ref{AcyclicSmallN}, we will show that $\mathcal{C}^{(g, n)}_*$ is acyclic when $n \leq 1$.
Denote by $\mathcal{R}^{(g, n)}_*$ the sub-chain complex of $\mathcal{C}^{(g, n)}_*$ spanned by those pairs $(\mathbf{P} \to \mathbf{C}, \omega)$ such that $\mathbf{P}$ has a vertex $v$ that has at least two markings from $\{1, \ldots, n\}$; this is the chain complex which computes the reduced rational homology of $\Theta_{g, n}^{\mathrm{rep}}$. Denote by $\mathcal{Q}^{(g, n)}_*$ the sub-chain complex of $\mathcal{C}^{(g, n)}_*$ spanned by those pairs $(\mathbf{P} \to \mathbf{C}, \omega)$ where $\mathbf{C}$ has at least one vertex $v$ with weight at least $3$.

% We will prove the following result, which is equivalent to Theorems \ref{AcyclicSubcomplexes} and \ref{AcyclicSmallN}.
% \begin{theorem}\label{AcyclicSummary}
%     The chain complexes $\mathcal{Q}^{(g, n)}_*$, $\mathcal{R}^{(g, n)}_*$ and $\mathcal{Q}^{(g, n)}_* \cap \mathcal{R}^{(g, n)}_*$ are acyclic for all $g \geq 2$ and all $n$. The chain complex $\mathcal{C}^{(g, n)}_*$ is acyclic for all $g \geq 2$ and $n \leq 1$.
% \end{theorem}

We will show that 
the chain complexes $\mathcal{R}^{(g, n)}_*$ and $\mathcal{Q}^{(g, n)}_* \cap \mathcal{R}^{(g, n)}_*$ are acyclic for all $g \geq 2$ and all $n \geq 2$ (Theorem \ref{RepAcyclic}),
  that the chain complex $\mathcal{Q}^{(g, n)}_*$ is acyclic for all $g \geq 2$ and all $n \geq 0$ (Theorem \ref{3endAcyclic}),  
   and that the chain complex $\mathcal{C}^{(g, n)}_*$ is acyclic for all $g \geq 2$ and $n \leq 1$ (Theorem \ref{OneMarkingAcyclic}).
Thus, Theorem \ref{RepAcyclic} and Theorem \ref{3endAcyclic} prove Theorem \ref{AcyclicSubcomplexes}, and Theorem \ref{OneMarkingAcyclic} gives part (1) of Proposition \ref{mainthm:smalln}.

The proofs of these theorems are informed by previous work of Chan--Galatius--Payne on contractibility criteria for symmetric $\Delta$-complexes \cite{cgp-marked}, as well as work of Conant--Gerlits--Vogtmann \cite{conant-gerlits-vogtmann-cut} on the acyclicity of the subcomplex of Kontsevich's graph complex spanned by graphs with cut vertices.

\subsection{The homology of $\Theta_{g, n}^{\mathrm{rep}}$}
It will be useful to isolate specific types of edges of covers with repeated markings.

\begin{figure}[h]
    \centering
    \includegraphics[scale=1.25]{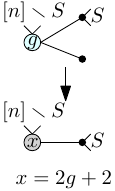}
    \caption{The cover $\mathbf{B}_S \to \mathbf{E}_S$.}
    \label{supportingedgeS}
\end{figure}

\begin{definition}
    For a $\Gamma^{\cH}_{g, n}$-object $\mathbf{P} \to \mathbf{C}$ with repeated markings, we say an edge $e \in E(\mathbf{C})$ is a \textit{supporting edge}, with support equal to $S \subseteq [n]$, if, upon contracting all edges of $\mathbf{C}$ which are not equal to $e$, as well as their preimages in $\mathbf{P}$, we obtain the cover $\mathbf{B}_S \to \mathbf{E}_S$ depicted in Figure \ref{supportingedgeS}. If $|S| = i$, we will call $e$ an \textit{$i$-supporting edge}.  Note that necessarily $i\ge 2$.
\end{definition} 

\begin{definition}
Given a $\Gamma_{g, n}^{\cH}$-object  $\mathbf{P} \to \mathbf{C}$, we define the \textit{supporting edge retraction} of $\mathbf{P} \to \mathbf{C}$ to be the cover obtained by contracting all supporting edges in $\mathbf{C}$ and their preimages in $\mathbf{P}$.
\end{definition}
\begin{theorem}\label{RepAcyclic}
    For all $g\geq 2$ and $n \geq 2$, the chain complexes $\mathcal{R}^{(g, n)}_*$ and $\mathcal{R}^{(g, n)}_* \cap \mathcal{Q}^{(g, n)}_*$ are acyclic.
\end{theorem}
\begin{proof}
We will prove the theorem only for $\cR^{(g, n)}_*$, as the same argument works for $\mathcal{R}^{(g, n)}_* \cap \mathcal{Q}^{(g, n)}_*$. For ease of notation, fix $g, n \ge 2$ and put \[\cR_* := \cR^{(g,n)}_*. \] First, filter $\cR_*$ as follows: let
\[ \cR^{\geq i}_* \hookrightarrow \cR_* \]
be the subcomplex generated by covers which have a $k$-supporting edge for some $k \geq i$. More precisely, we mean that $\cR^{\geq i}_*$ is spanned by covers obtained by edge-contraction from covers with supporting edges of this type. We apply this definition even when $i=n+1$, in which case $\cR^{\ge {n+1}}_* = 0$.  Then we have a filtration
\[ 0 = \cR^{\ge n+1}_* \hookrightarrow  \cR^{\geq n}_* %\hookrightarrow \cdots \hookrightarrow \cR_*^{\geq i + 1} \hookrightarrow \cR_*^{\geq i} \hookrightarrow \cR_*^{\geq i - 1} 
\hookrightarrow \cdots \hookrightarrow \cR^{\geq 2}_* = \cR_*.  \]
Passing to the associated spectral sequence, it suffices to show that for each $i=2,\ldots,n$, the successive quotient chain complexes
\[ \cR^{i}_* := \cR^{\geq i}_* / \cR^{\geq i + 1}_* \]
are acyclic. These quotient chain complexes are spanned by covers with $i$-supporting edges and their edge-contractions, but do not include any covers with $k$-supporting edges or their edge contractions for any $k > i$. Now we filter $\cR^i_*$. Define
\[F_p\cR^i_* \hookrightarrow \cR^i_*  \]
to be the sub-chain complex spanned by graphs with at most $p$ non-supporting edges.  
The number of non-supporting edges cannot increase under edge contraction, so $F_p \cR^i_*$ really is a subcomplex.  We obtain
 an ascending filtration
\[0 = F_{-1}\cR^i_* \hookrightarrow F_{0}\cR^i_* \hookrightarrow \cdots\hookrightarrow \cR^i_*\]
and again by considering the associated spectral sequence, it suffices to show that successive quotients
\[ G_p \cR^i_* := F_p \cR^i_* / F_{p -1}\cR^i_*  \]
are acyclic, in order to conclude that $\cR^i_*$ and hence $\cR_*$ is acyclic. For fixed $i$ and $p$, let $A_{i, p}$ denote the set of isomorphism classes of $\Gamma^{\cH}_{g, n}$-objects $\P \to \C$ where $|E(\C)| = p$ and which (1) do not have any supporting edges, (2) admit a contraction from a cover with an $i$-supporting edge, and (3) do not admit a contraction from any covers with $k$-supporting edges for $k > i$. Then we have a direct sum decomposition
\[ G_p \cR^i_* = \bigoplus_{\P \to \C \in A_{i,p}} \cL^{\P \to \C}_*, \]
where $\cL^{\P \to \C}_*$ is the sub-chain complex consisting of those covers whose supporting edge retraction is equal to $\P \to \C$. This direct sum decomposition holds because the differential on $G_p \cR^i_*$ is given by a signed sum of supporting edge contractions, and hence preserves the supporting edge retraction of a given cover. Next, given $\P \to \C \in A_{i, p}$, we have a tensor product decomposition
\[\cL_*^{\P \to \C} \cong \left(\bigotimes_{v \in V^{\mathrm{rep}}_i(\P)} (\QQ \xrightarrow{\sim} \QQ)\right)[1 - p], \]
where $V^{\mathrm{rep}}_i(\P)$ denotes the set of vertices of $\P$ which contain exactly $i$ markings, and the first copy of $\QQ$ is in degree 1. This tensor product decomposition holds because a generator of $\cL_*^{\P \to \C}$ is determined by a choice of subset of those vertices of $\P$ which contain $i$ markings: the corresponding generator is determined by expanding a single $i$-supporting edge from the image of each chosen vertex in $\C$.  (Since $i\ge 2$, such an expansion is indeed possible, producing a {\em stable} $S$-marked target tree.)  The degree shift is required to account for the $p$ edges of $\P \to \C$.
Altogether, this shows that $\cL^{\P \to \C}_*$ is a tensor product of acyclic chain complexes, so $\cL^{\P \to \C}_*$ is itself acyclic, and the proof is complete.
\end{proof}

\subsection{The homology of $\Theta_{g, n}^{\geq 3}$} 
\begin{figure}[h]
    \centering
    \includegraphics[scale=1.25]{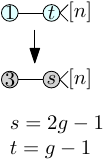}
    \caption{The cover $\mathbf{D} \to \mathbf{F}$.}
    \label{3end}
\end{figure}

We will now show that the chain complex $\mathcal{Q}^{(g, n)}_*$ is acyclic. It will again be convenient to name particular types of edges.
\begin{definition}
Suppose $\mathbf{P} \to \mathbf{C}$ is an object of $\Gamma^{\mathcal{H}}_{g, n}$, and that $\mathbf{C}$ has a vertex of weight at least $3$. 
\begin{enumerate}
\item We say $e \in E(\mathbf{C})$ is a \textit{$3$-end} if upon contracting all edges in $\mathbf{C}$ except for $e$, and their preimages in $\mathbf{P}$, we obtain the cover $\mathbf{D} \to \mathbf{F}$ in Figure \ref{3end}.
\item We say a cover $\mathbf{P}' \to \mathbf{C}'$ is a \textit{3-end expansion} of $\mathbf{P} \to \mathbf{C}$ if $\mathbf{P} \to \mathbf{C}$ is obtained from $\mathbf{P}' \to \mathbf{C}'$ by contracting a sequence of $3$-ends.
\end{enumerate}
\end{definition}

\begin{figure}[h]
\centering
\includegraphics[scale=1.25]{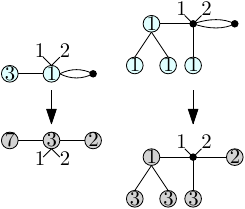}
\caption{A cover in $\Gamma^{\mathcal{H}}_{5, 2}$ and its maximal expansion by $3$-ends.}
\label{Weight3ExpansionFig}
\end{figure}

 It is straightforward to see that for any cover $\mathbf{P} \to \mathbf{C}$, the poset of $3$-end expansions of $\mathbf{P} \to \mathbf{C}$ has a maximal element, as in the following lemma. We omit the %straightforward 
 %% --we just said straightforward --MC
 proof: see Figure \ref{Weight3ExpansionFig} for an example of how this expansion is constructed.
 \begin{lemma}\label{3EndExpansion}
     Let $\mathbf{P} \to \mathbf{C}$ be an object of $\Gamma^{\cH}_{g, n}$. Then the poset of $3$-end expansions of $\mathbf{P} \to \mathbf{C}$ has a unique maximal element $\mathbf{P}' \to \mathbf{C}'$, and this expansion is canonical in the sense that any automorphism of $\mathbf{P} \to \mathbf{C}$ lifts to an automorphism of $\mathbf{P}' \to \mathbf{C}'$.
 \end{lemma}

Given a $\Gamma^{\cH}_{g, n}$-object $\mathbf{P} \to \mathbf{C}$, let $A(\mathbf{P} \to \mathbf{C})$ be the set of isomorphism classes of covers obtained from $\mathbf{P} \to \mathbf{C}$ by contracting $3$-ends. We define a chain complex $\mathcal{Q}^{\mathbf{P} \to \mathbf{C}}_*$ as follows: the vector space $\mathcal{Q}^{\mathbf{P} \to \mathbf{C}}_p$ is spanned by pairs $(\mathbf{H} \to \mathbf{K}, \omega)$, where $\mathbf{H} \to \mathbf{K}$ is an element of $A(\mathbf{P} \to \mathbf{C})$ with $|E(\mathbf{K})| = p + 1$, and $\omega \colon  [p] \to E(\mathbf{K})$ is an edge-labelling. These generators are subject to the usual relation
\[(\mathbf{H} \to \mathbf{K}, \omega \circ \rho ) = \mathrm{sgn}(\rho) (\mathbf{H} \to \mathbf{K}, \omega)\]
for $\rho \in \Aut([p])$. The differential on $\mathcal{Q}^{\mathbf{P} \to \mathbf{C}}_*$ is given by the signed sum of $3$-end contractions; we set it equal to $0$ on any generators which do not have any $3$-ends.

\begin{proposition}\label{Local3end}
Suppose $\mathbf{P} \to \mathbf{C}$ has a $3$-end and is maximal with respect to expanding $3$-ends. Then $\mathcal{Q}^{\mathbf{P} \to \mathbf{C}}_*$ is acyclic.
\end{proposition}

\begin{proof}
    First consider the case where $\mathbf{C}$ has no automorphisms. This implies that all $3$-end contractions of $\mathbf{C}$ have no automorphisms, since any automorphism of the target tree of a $\Gamma^{\cH}_{g, n}$-object must lift to an automorphism of its maximal $3$-end expansion. Let $q + 1$ be the number of distinct $3$-ends of $\mathbf{C}$. We can understand $\mathcal{Q}^{\mathbf{P} \to \mathbf{C}}_*$ as a shift of the augmented cellular chain complex of the standard $q$-simplex $\sigma^{q}$, viewed as the space parameterizing assignments of nonnegative lengths to the $q + 1$ distinct $3$-ends of $\mathbf{C}$, such that the lengths sum to one. Note that $q\ge 0$ by assumption, so that $\sigma^q$ is nonempty.  So in the automorphism-free case, $\mathcal{Q}^{\mathbf{P} \to \mathbf{C}}_*$ is acyclic. 

    For the general case, when $\mathbf{C}$ and its contractions may have automorphisms, fix a labelling of the edges of $\mathbf{C}$, and denote the resulting object by $\mathbf{C}^\dagger$. This induces a labelling of the edges of each contraction of $\mathbf{C}$. Let $A( \mathbf{C}^\dagger)$ be the set consisting of $\mathbf{C}^\dagger$ and all of its $3$-end contractions. We can make a chain complex $\mathcal{Q}^{\mathbf{C}, \dagger}_{*}$ which in degree $p$ is spanned by pairs $[\mathbf{K}, \omega]$ where $\mathbf{K}$ is an element of $A(\mathbf{C}^{\dagger})$ with $|\mathbf{K}| = p + 1$, and $\omega\colon  [p] \to E(\mathbf{K})$ is a bijection, subject to the usual relations under the action of $\Aut([p])$. Observe that there is a canonical action of $\Aut(\mathbf{C})$ on the chain complex $\mathcal{Q}^{\mathbf{C}, \dagger}_{*}$, and $\mathcal{Q}^{\mathbf{P} \to \mathbf{C}}_*$ is identified with the $\Aut(\mathbf{C})$-coinvariants of the complex $\mathcal{Q}^{\mathbf{C}, \dagger}_{*}$, by the second part of Lemma \ref{3EndExpansion}. Since $\Aut(\mathbf{C})$ is finite, it has no homology over the rationals. Moreover, $\mathcal{Q}^{\mathbf{C}, \dagger}_*$ is acyclic by the first part of the proof. We conclude that %\[H_*\left(\left(\mathcal{Q}^{\mathbf{C}, \dagger}_{*}\right)_{\Aut(\mathbf{C})}\right) = \left(H_*\left(\mathcal{Q}^{\mathbf{C}, \dagger}_{*}\right)\right)_{\Aut(\mathbf{C})} = 0, \]
    \[H_*(\mathcal{Q}^{\mathbf{P} \to \mathbf{C}}_*) = H_*((\mathcal{Q}^{\mathbf{C}, \dagger}_{*})_{\Aut(\mathbf{C})}) = (H_*(\mathcal{Q}^{\mathbf{C}, \dagger}_{*}))_{\Aut(\mathbf{C})} = 0, \]
    as desired.
    \end{proof}
    %Since no graphs in $A(\bfT^{\dagger})$ have automorphisms, we see that $\mathcal{Q}^{\bfT, \dagger}_{*}$ is acyclic. There is an action of $\Aut(\bfT)$ on $\mathcal{Q}^{\bfT, \dagger}_{*}$, induced by the action of $\Aut(\bfT)$ on $E(\bfT)$: given $\psi \in \Aut(\bfT)$ and an edge-labelling $\omega: E(\bfT^\dagger) \to [r]$, we set $\psi \cdot (\bfT^\dagger, \omega) = (\bfT^\dagger, \omega \circ \psi)$. For any $\mathbf{P} \in A(\bfT^\dagger)$, there is a unique set of edges $e_1, \ldots, e_\ell \in E(\bfT^\dagger)$ such that $\mathbf{P}$ was obtained from $\bfT^\dagger$ by contracting $e_1, \ldots, e_\ell$. We set $\psi \cdot (\mathbf{P}, \omega) = (\bfT^{\dagger}/\{\psi(e_1, \ldots, \psi(e_\ell) \}, \omega \circ \kappa)$, where $\kappa$ is the unique bijection making the diagram
    % \[ \begin{tikzcd}
    %     &E(\bfT^{\dagger})  &E(\bfT^{\dagger})\arrow[l, "\psi"] \\
    %     &E(\mathbf{P})\arrow[u]   &E(\bfT^\dagger / \{\psi(e_1), \ldots, \psi(e_k)\})\arrow[l, "\kappa"] \arrow[u]
    % \end{tikzcd} \]
    % commute. 

We now prove that $\mathcal{Q}^{(g, n)}_*$ is acyclic.

\begin{theorem}\label{3endAcyclic}
    For $g \geq 2$ and $n \geq 0$, the chain complex $\mathcal{Q}^{(g, n)}_*$ is acyclic.
\end{theorem}
\begin{proof}
    Let $F_p \mathcal{Q}^{(g, n)}_*$ denote the subspace spanned by those covers whose target tree has at most $p$ edges which are not $3$-ends. This defines a bounded, increasing filtration of $\mathcal{Q}^{(g, n)}_*$. The $E^0$ page
    \[E^{0}_{p, q} = F_p \mathcal{Q}^{(g, n)}_{p + q}/ F_{p - 1}\mathcal{Q}^{(g, n)}_{p + q} \]
    of the associated spectral sequence is spanned by covers whose target tree has exactly $p$ edges which are not $3$-ends. The differential $\partial_0\colon E^{0}_{p, q} \to E^{0}_{p, q-1}$ is given by a signed sum of $3$-end contractions. Therefore, by Lemma \ref{3EndExpansion}, the $p$th column of the $E^0$ page breaks up into a direct sum of chain complexes of the form $\mathcal{Q}^{\mathbf{P} \to \mathbf{C}}_*$, where $\mathbf{C}$ has at least one $3$-end, and the tree obtained from $\mathbf{C}$ by contracting all $3$-ends has $p$ edges. Proposition \ref{Local3end} then implies that the $E^1$ page vanishes, which completes the proof.
    \end{proof}

\subsection{Calculations on $\Theta_{g, n}$ for $n \leq 2$} 
\begin{figure}[h]
        \centering
        \includegraphics[scale=1.25]{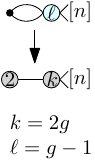}
        \caption{The cover $\mathbf{J} \to \mathbf{K}$.}
        \label{2end}
\end{figure}

We conclude this section by proving Proposition \ref{mainthm:smalln}. The first part of Proposition \ref{mainthm:smalln} asserts that $\mathcal{C}^{(g, n)}_*$ is acyclic for $n \leq 1$, and the proof is similar to the one that $\mathcal{Q}^{(g, n)}_*$ is acyclic. Once again, we isolate particular types of edges:
\begin{definition}
Let $\mathbf{P} \to \mathbf{C}$ be a $\Gamma_{g, n}^{\cH}$-object. An edge $e \in E(\mathbf{C})$ is called a \textit{$2$-end} if upon contracting all edges of $\mathbf{C}$ except for $e$, and their preimages in $\mathbf{P}$, we obtain the cover $\mathbf{J} \to \mathbf{K}$ in Figure \ref{2end}.
\end{definition}
The key to the proof of acyclicity of $\mathcal{C}^{(g, n)}_*$ when $n \leq 1$ is the following lemma.

\begin{figure}[h]
    \centering
\includegraphics[scale=1.25]{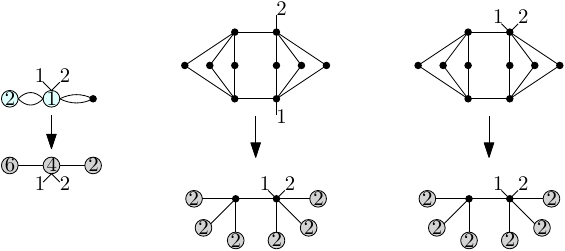}
    \caption{A cover $\mathbf{P} \to \mathbf{C}$ in $\Gamma_{5, 2}^{\cH}$ and the two distinct maximal elements of its poset of $2$-end expansions. When $n \leq 1$, this poset always has a unique maximal element, as explained in the proof of Lemma \ref{2EndExpansions}.}
    \label{noncanonical_expansions}
\end{figure}

\begin{lemma}\label{2EndExpansions}
    Let $\mathbf{P} \to \mathbf{C}$ be an object of $\Gamma^{\cH}_{g, n}$ for $n \leq 1$. Then the poset of expansions of $\mathbf{P} \to \mathbf{C}$ by $2$-ends has a unique maximal element $\mathbf{P}' \to \mathbf{C}'$. Moreover, this expansion is canonical in the sense that any automorphism of $\mathbf{P} \to \mathbf{C}$ lifts to one of $\mathbf{P}' \to \mathbf{C}'$.
\end{lemma}
\begin{proof}
It is clear how to construct the graph $\mathbf{C}'$: for every vertex of $\mathbf{C}$ with weight $d \geq 2$, one expands $\lfloor d/2 \rfloor$ many $2$-ends from $v$, leaving behind a vertex of weight $d - 2\lfloor d/2 \rfloor$ (if $d = 2$, this expansion should only be performed if it preserves the stability condition, i.e. only if the vertex is not already part of a $2$-end). This uniquely determines a cover $P'$, but does not determine the marking function on $P'$. If $n = 0$, then there is no marking function, so $\mathbf{P}'$ is determined. For $n = 1$, the only ambiguity arises when $v$ supports the unique marking, and the preimage of $v$ in $\mathbf{C}'$ has $2$ preimages in the graph $P'$, so one has to make a choice as to which fiber to mark. However, since $n = 1$, both choices are equivalent, as they differ by the $\ZZ/2\ZZ$-action on $P'$. Therefore $\mathbf{P}'$ is also determined when $n = 1$. The statement on lifting of automorphisms is straightforward to check. The lemma fails when $n > 1$, because in general there is no canonical way of distributing the markings supported at $v$ among the fibers over $v$ in $P'$. See Figure \ref{noncanonical_expansions} for an example.
\end{proof}

\begin{figure}
    \centering
    \includegraphics[scale=1]{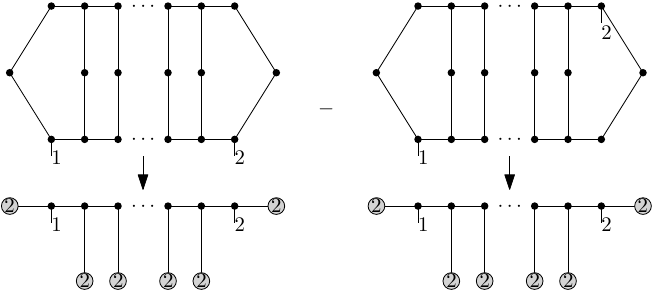}
    \caption{A cycle spanning $\widetilde{H}_{2g}(\Theta_{g, 2} ; \QQ)$}
    \label{fig:H_g2-cycle}
\end{figure}

Given Lemma \ref{2EndExpansions}, the proof of the following theorem is completely analogous to the proof of Theorem \ref{3endAcyclic}; we will only outline the necessary steps.
\begin{theorem}\label{OneMarkingAcyclic}
    For $g\geq 2$ and $n \leq 1$, the chain complex $\mathcal{C}^{(g, n)}_*$ is acyclic.
\end{theorem}
\begin{proof}
    First, define $B(\mathbf{P} \to \mathbf{C})$ to be the set of isomorphism classes of $\Gamma_{g, n}^{\cH}$-objects obtained from $\mathbf{P} \to \mathbf{C}$ by contracting $2$-ends. Then, use this to define a chain complex $\cG^{\mathbf{P} \to \mathbf{C}}_*$ analogously to $\cQ^{\mathbf{P} \to \mathbf{C}}_*$, where the differential is given by a signed sum of $2$-end contractions. The proof that $\cG^{\mathbf{P} \to \mathbf{C}}_*$ is acyclic, for $\mathbf{P} \to \mathbf{C}$ maximal with respect to expanding $2$-ends, is exactly the same as the proof of Proposition \ref{Local3end}. Finally, one proves the theorem by filtering $\mathcal{C}^{(g, n)}_*$: set $F_p \mathcal{C}^{(g, n)}_*$ to be the subcomplex of $\mathcal{C}^{(g, n)}_*$ spanned by those covers with at most $p$ edges which are not $2$-ends. Then the $p$th column of the $E^0$ page of the associated spectral sequence breaks up into a direct sum of complexes of the form $\cG^{\mathbf{P} \to \mathbf{C}}_*$ by Lemma \ref{2EndExpansions}, so the $E^1$ page vanishes, and the result follows.
\end{proof}
Theorem \ref{OneMarkingAcyclic} gives part (1) of Proposition \ref{mainthm:smalln}. Part (2) states that
\[W_0 H^{2g + 1}_c(\cH_{g, 2};\QQ) \cong \QQ, \]
and that the corresponding $S_2$-representation is trivial if $g$ is even, and given by the sign representation if $g$ is odd. We prove this now by writing down an explicit cycle in $\cC^{(g, 2)}_{2g}$ corresponding to this class. See Figure \ref{fig:H_g2-cycle}.

\begin{proof}[Proof of Proposition \ref{mainthm:smalln}, part (2)]
We have an isomorphism of $S_2$-representations
\[W_0 H^{2g + 1}_c(\cH_{g, 2};\QQ) \cong \widetilde{H}_{2g}(\Theta_{g, 2};\QQ)^\vee\]
by Corollary \ref{cor:WeightZeroIdentification}. We have
\[\widetilde{H}_{2g}(\Theta_{g, 2};\QQ) = H_{2g}\left(\cC^{(g, 2)}_*\right). \]
Observe that $2g$ is the top homological degree of $\cC^{(g, 2)}$: the maximal number of edges of a stable tree with $2g + 4$ legs is $2g + 1$. Therefore, any cycle in $\cC^{(g, 2)}_{2g}$ defines a class in homology. Any target tree for a cover in $\cC^{(g, 2)}_{2g}$ must be trivalent, and to be a nonzero element, it cannot have any automorphisms which act by an odd permutation of the edge set. It is straightforward to conclude that such a tree must be equal to the tree depicted in Figure \ref{fig:H_g2-cycle}. This tree $\C$ has two covers, depicted in Figure \ref{fig:H_g2-cycle}. Therefore $\dim \cC^{(g, 2)}_{2g} = 2$, where a basis is given by choosing any edge-labelling of the aforementioned tree. One can verify directly that neither one of these basis elements form a cycle on their own, but their difference does. From this we conclude that $H_{2g}(\cC_*^{(g,2)}) \cong \QQ$. To understand the $S_2$-representation, we note that when $g$ is even, the transposition in $S_2$ induces an even permutation on any edge-labelling of the given tree, and when $g$ is odd, the transposition induces an odd permutation of the edge labels. 
\end{proof}

\begin{remark}\label{remark: generalizedrepeatedmarkings}
Theorem \ref{RepAcyclic} generalizes  to other spaces of admissible covers. Fix an integer $N > 0$ and let $G$ be a prime cyclic group, which we now write additively to be consistent with our notation for $G=\ZZ/2\ZZ$. Let
\[\mu\colon  \{w_1, \ldots, w_N\} \to G  \]
 be a function such that the image of $\mu$ generates $G$, which additionally satisfies
 \[ \sum_{i = 1}^{N} \mu(w_i) = 0 \]
 where $0 \in G$ denotes the identity element. For any integer $n \geq 0$, we can extend $\mu$ to a function
 \[ \{1, \ldots, n\} \cup \{w_1, \ldots, w_N\} \to G \]
 by setting the image of each $i \in \{1, \ldots, n\}$ to be $0$; for ease of notation, we will also call this extension $\mu$. We set the notation
\[\overline{\cM}_{0, n + N}^G(\mu):= \overline{\cM}_{0, \{1, \ldots, n \} \cup \{w_1, \ldots, w_N \}}^{G}(\mu) \]
and define $\cM_{0, n + N}^G(\mu)$ similarly. We now define an intermediate locus
\[\cM_{0, n + N}^G(\mu) \subset \widetilde{\cM}_{0, n + N}^G(\mu) \subset \overline{\cM}_{0, n + N}^G(\mu) \]
in analogy with the space $\widetilde{\cH}_{g, n}$ of $n$-marked hyperelliptic curves of genus $g$ together with a labelling of their Weierstrass points, considered in \S\ref{section:tildeH_gn}. Given a graph-theoretic pointed admissible $G$-cover $\P \to \C \in \mathrm{Ob}(\Gamma_{0, n+N}^G(\mu))$, where $\Gamma_{0, n+N}^{G}(\mu)$ is the category defined in Definition \ref{def:Gamma-category}, we say that $\P \to \C$ is \textit{forbidden} if all of the following conditions hold:
\begin{enumerate}[(a)]
\item $|E(\C)| = 1$,
\item If we erase all of the legs labelled by $\{1, \ldots, n\}$ from $\C$, the resulting $\{w_1, \ldots, w_N\}$-marked tree is not stable in the sense of Definition \ref{def:GraphicCovers}, and
\item the source graph $\P$ has no vertices supporting repeated markings among $\{1, \ldots, n\}$.
\end{enumerate}
Each forbidden cover $\P \to \C$ corresponds uniquely to a boundary divisor of $\overline{\cM}_{0, n + N}^G(\mu)$, and we define $\widetilde{\cM}_{0, n + N}^G(\mu)$ to be the complement in $\overline{\cM}_{0, n + N}^G(\mu)$ of those boundary divisors which are not forbidden.

When $G = \ZZ/2\ZZ = \{0, 1\}$, $N = 2g + 2$, and $\mu(w_i) = 1$ for all $i$, the forbidden divisors are precisely those of type (1) and (2) in Definition \ref{defn:orderedH_gn}, and we have
\[ \widetilde{\cM}_{0, n + 2g + 2}^{\ZZ/2\ZZ}(\mu) \cong \wt{\cH}_{g, n}. \]
For general $G$ and $\mu$, the space $\widetilde{\cM}_{0, n + N}^G(\mu)$ can be identified with the moduli space of smooth $N$-pointed admissible $G$-covers of $\PP^1$, with monodromy specified by $\mu$, together with $n$ distinct marked points on the source curve. This space admits an $S_n$-action given by permuting the $n$ marked points on the source, and the isomorphism with $\widetilde{\cM}_{0, n + N}^G(\mu)$ is $S_n$-equivariant.

The dual complex $\widetilde{\Theta}_{0, n+N}^G(\mu)$ of the normal crossings compactification
\[ \widetilde{\cM}_{0, n + N}^G(\mu) \subset \overline{\cM}_{0, n + N}^G(\mu) \]
is the subcomplex of $\Delta_{0, n+ N}^G(\mu)$ of those simplices which have no forbidden vertices. The analogue of Theorem \ref{RepAcyclic} holds for $\widetilde{\Theta}_{0, n+N}^G(\mu)$: the subcomplex parameterizing graph-theoretic admissible $G$-covers $\P \to \C$ where $\P$ has a repeated marking is acyclic. Our proof of Theorem \ref{RepAcyclic} carries through to this setting, mutatis mutandis. In Remark \ref{remark: generalizedgraphsum} below, we explain how this leads to a generalization of Theorem \ref{mainthm:FrobChar} for these spaces.
\end{remark}

\section{A graph sum formula for $\mathsf{h}_g$}\label{section: graphsumformula}

Recall from the introduction that
\[ \mathsf{h}_g = \sum_{n \geq 0} \sum_{i  = 0}^{4g - 2 + 2n} (-1)^i \ch_n W_0 H^i_c (\cH_{g, n}; \QQ) \in \hat{\Lambda} \]
denotes the generating function for the weight zero equivariant Euler characteristics of the moduli spaces $\cH_{g, n}$. In this section we will prove Theorem \ref{mainthm:FrobChar}, thus establishing our sum-over-graphs formula for $\mathsf{h}_g$. We let $T_{2g + 2}$ denote the set of isomorphism classes of stable trees with $2g + 2$ unlabelled legs. When each leg is given monodromy marking equal to $1 \in \ZZ/2\ZZ$, such a tree $C$ has a unique graph-theoretic admissible $\ZZ/2\ZZ$-cover $P_C\to C$. Let $T_{2g + 2}^{<3}$ denote the subset of $T_{2g + 2}$ consisting of those trees such that no vertex supports more than two leaves, and for a tree $C$ we write $E_C$ for its set of edges. We restate Theorem \ref{mainthm:FrobChar} for convenience.
\begin{customthm}{A}
    We have
    \[ \mathsf{h}_g = \sum_{C \in T_{2g + 2}^{<3}} %\mathsf{h}_{P_C, E_C, \mathrm{Aut}(P_C)}, 
      \frac{(-1)^{|E_C|}}{|\Aut(P_C)|} \sum_{\tau \in \Aut(P_C)} \mathrm{sgn}(\tau|_{E_C}) \prod_{k \geq 1} (1 + p_k)^{f(P_C, \tau, k)}
    \]
    where $E_C$ is the set of edges of the tree $C$, 
    $p_k = \sum_{n > 0} x_n^{k} \in \hat{\Lambda}$ is the $k$th power sum symmetric function, and $k \cdot f(P_C, \tau, k)$ is the compactly supported Euler characteristic of the set of points in $P_C$ which have orbit of length $k$, under the action of $\tau$.
\end{customthm}

We will prove Theorem \ref{mainthm:FrobChar} through a series of intermediate results. Throughout this section, we tacitly replace the symmetric $\Delta$-complex $\Theta_{g, n}$ with its geometric realization. 

\begin{lemma}
We have
\[ \mathsf{h}_g = - \sum_{n \geq 0} \chi_c^{S_n}(\Theta_{g, n} \smallsetminus (\Theta_{g, n}^{\mathrm{rep}} \cup \Theta_{g, n}^{\geq 3})), \]
where $\chi_c^{S_n}( \cdot )$ denotes the $S_n$-equivariant compactly supported Euler characteristic.
\end{lemma}
\begin{proof}
Via the identification 
\[W_0 H^i_c(\cH_{g, n} ; \QQ) \cong \wt{H}_{i -1}(\Theta_{g, n}; \QQ)^\vee \]
% \maddie{we spend a lot of time constructing a normal crossings compactification, and then I don't think we ever say that is why this identification holds / there is no formal statement of this identification. I don't know if we need to change anything, but one suggestion is that we could add a statement of this to the end of section 4? this would also motivate section 5: why are we now interested in the homology of theta gn? also i think this identification is used in a few places in section 5 as well.}
of Corollary \ref{cor:WeightZeroIdentification}, and using that the Frobenius characteristic of a representation of $S_n$ equals that of its dual, we can write
\begin{align*}
    \mathsf{h}_g &= \sum_{n \geq 0} \sum_{i  = 0}^{4g - 2 + 2n} (-1)^i \ch_n \wt{H}_{i - 1} (\Theta_{g, n}; \QQ) \\&= \sum_{n \geq 0} -\wt{\chi}^{S_n}(\Theta_{g, n}),
\end{align*}
where $\wt{\chi}^{S_n}(\cdot)$ denotes the $S_n$-equivariant reduced Euler characteristic. Since $\Theta_{g, n}$ is connected and compact, and $S_n$ acts trivially on $H_0(\Theta_{g, n};\QQ) \cong \QQ$, we have
\begin{align*}
   -\sum_{n \geq 0} \wt{\chi}^{S_n}(\Theta_{g, n}) &= \sum_{n \geq 0} h_n - \sum_{n \geq 0}\chi_c^{S_n}(\Theta_{g, n}),
\end{align*}
where $h_n \in \Lambda$ is the $n$th homogeneous symmetric function, defined as the Frobenius characteristic of the trivial $S_n$-representation. By the additivity of the compactly supported Euler characteristic under stratification, we can write
\begin{align*}
   \sum_{n \geq 0}\chi_c^{S_n}(\Theta_{g, n}) &= \sum_{n \geq 0} \left( \chi_c^{S_n}(\Theta_{g, n} \smallsetminus (\Theta_{g, n}^{\mathrm{rep}} \cup \Theta_{g, n}^{\geq 3})) + \chi_c^{S_n}(\Theta_{g, n}^{\mathrm{rep}} \cup \Theta_{g, n}^{\geq 3})\right)
\end{align*}   
Since the union $\Theta_{g, n}^{\mathrm{rep}} \cup \Theta_{g, n}^{\geq 3}$ is compact and connected, with vanishing reduced rational homology by Theorem \ref{AcyclicSubcomplexes}, and $S_n$ acts trivially on $H_0(\Theta_{g, n}^{\mathrm{rep}} \cup \Theta_{g, n}^{\geq 3} ; \QQ)$, we have
\[  \chi_c^{S_n}(\Theta_{g, n}^{\mathrm{rep}} \cup \Theta_{g, n}^{\geq 3}) = h_n, \]
and the proof is complete.
\end{proof}

\begin{lemma}\label{hgStratificationlemma}
We have
\[ \mathsf{h}_g = -\sum_{C \in T_{2g + 2}^{<3}} \sum_{n \geq 0} \chi_c^{S_n}\left( \left(\mathrm{Conf}_n(P_C) \times (\Delta^{|E_C| - 1})^\circ \right)/\Aut(P_C)\right). \]
\end{lemma}
\begin{proof}
    We can stratify the space
\[X_{g, n} : = \Theta_{g, n} \smallsetminus (\Theta_{g, n}^{\mathrm{rep}} \cup \Theta_{g, n}^{\geq 3})  \]
by the $\Gamma_{g}^{\cH}$-object that arises when we forget the marking function and delete all legs with monodromy equal to $0$, as well as their preimages, and then stabilizing: stabilization process only entails the removal of 2-valent vertices, because we are outside the locus $\Theta_{g, n}^{\mathrm{rep}}$. Such an object is uniquely specified by an element $C$ of $T_{2g + 2}^{<3}$, which determines its covering $P_C$. Since we have removed the repeated marking locus, the stratum corresponding to $P_C \to C$  is $S_n$-equivariantly homeomorphic to
\[ \left(\mathrm{Conf}_n(P_C) \times (\Delta^{|E_C| - 1})^\circ \right)/\Aut(P_C).  \]
Above, $(\Delta^{|E_C| - 1})^\circ$ denotes the interior of the standard $|E_C| - 1$ simplex $\Delta^{|E_C| - 1}$, viewed as the space parameterizing metrics $\ell\colon  E_C \to \RR_{> 0}$ of total length one. The space $\mathrm{Conf}_n(P_C)$ is the configuration space of $n$ distinct points on $P_C$, and the action of $\Aut(P_C)$ is diagonal: one finds that the morphism $P_C \to C$ can be reconstructed from $P_C$, so $\Aut(P_C)$ naturally acts on $C$ and hence on $|E_C|$ and $\left(\Delta^{|E_C|- 1}\right)^{\circ}$. 
\end{proof}

We now show how to calculate the terms in the sum, following Gorsky's calculation of the $S_n$-equivariant Euler characteristic of $\mathrm{Conf}_n(X)/G$, where $X$ is an algebraic variety and $G$ is a finite subgroup of its automorphism group \cite{gorsky-equivariant}.
\begin{proposition}\label{CWContribution}
    Let $X$ be a finite CW complex, and let $E$ be a finite set. Set \[\Delta^{\circ} = \left \{\ell\colon  E \to \RR_{>0} \mid \sum_{e \in E}\ell(e) = 1 \right \}.\] Let $G$ be a finite group acting on both $X$ and $E$, and set
    \[\mathsf{h}_{X, E, G} := \sum_{n \geq 0} \chi_c^{S_n}\left( \left(\mathrm{Conf}_n(X) \times \Delta^\circ \right)/G\right).\]
    Then
    \[\mathsf{h}_{X, E, G} = -\frac{(-1)^{|E|}}{|G|} \sum_{g \in G} \mathrm{sgn}\left(g|_{E}\right) \prod_{k \geq 1} (1 + p_k)^{\chi_c(X_k(g))/k}, \]
    where $X_k(g)$ denotes the set of points of $X$ which have orbit of length $k$ under the action of $g$.
\end{proposition}
Before proving Proposition \ref{CWContribution}, we need two intermediate lemmas.
\begin{lemma}\label{ConfigurationSpaceChi}
    Suppose that $X$ is any finite CW complex. Then
    \[ f(t) := \sum_{n \geq 0} \chi_c(\mathrm{Conf}_n(X)) \frac{t^n}{n!} = (1 + t)^{\chi_c(X)}. \]
\end{lemma}
\begin{proof}
    We have the identity
    \[ \chi_c(X^n) = \sum_{k = 1}^{n} S(n, k) \chi_c(\mathrm{Conf}_k(X)), \]
    where $S(n, k)$, the Stirling number of the second kind, counts the number of partitions of $n$ with $k$ parts. It follows that
    \[g(t): =\sum_{n \geq 0} \chi_c(X^n) \frac{t^n}{n!} = e^{\chi_c(X)t} \]
    is the \textit{Stirling transform} of $f$, so that $f(t) = g(\log(1+t)) = (1 + t)^{\chi_c(X)}$, as claimed.
\end{proof}
\begin{lemma}\label{sigmafixedpoints}
For any group $H$ acting on a space $Y$,  denote by
\[ [Y]^{h} \]
the set of fixed points of $h \in H$ acting on $Y$. Then, for $X$, $E$, and $G$ as above, and $\sigma\in S_n$, we have
\[ \chi_c\left(\left[\left(\mathrm{Conf}_n(X) \times \Delta^\circ\right)/G\right]^{\sigma}\right) = - \frac{(-1)^{|E|}}{|G|} \sum_{g \in G} \mathrm{sgn}(g|_{E}) \cdot \chi_c\left([\mathrm{Conf}_n(X)]^{g^{-1}\sigma}\right). \]
\end{lemma}
\begin{proof}
Define
\[S = \{(g, \ell, y) \in G \times \Delta^\circ \times \mathrm{Conf}_n(X) \mid g \cdot (\ell, y) = \sigma \cdot (\ell, y) \}. \] Then we have a map
\[S \to \left[\left(\mathrm{Conf}_n(X) \times \Delta^\circ\right)/G \right]^{\sigma}, \]
which takes $(g, \ell, y)$ to $(y, \ell)$. The fibers of this map are all nonempty and have cardinality equal to $|G|$,
so
\[\chi_c\left(\left[\left(\mathrm{Conf}_n(X) \times \Delta^\circ\right)/G\right]^{\sigma}\right) = \frac{1}{|G|} \chi_c(S). \]

On the other hand, the projection $S \to G$ has fiber over $g \in G$ isomorphic to
\[ [\Delta^\circ]^{g} \times [\mathrm{Conf}_n(X)]^{g^{-1}\sigma}. \]
Therefore we have
\[\chi_c\left(\left[\left(\mathrm{Conf}_n(X) \times \Delta^\circ\right)/G\right]^{\sigma}\right) = \frac{1}{|G|} \sum_{g \in G} \chi_c([\Delta^\circ]^{g}) \cdot \chi_c\left([\mathrm{Conf}_n(X)]^{g^{-1}\sigma} \right). \]
The proof is finished upon noting that $[\Delta^\circ]^g$ is again an open simplex, whose dimension modulo $2$ is equal to $|E| + \mathrm{sgn}(g|_{E}) - 1$.
\end{proof}

We can now prove Proposition \ref{CWContribution}.
\begin{proof}[Proof of Proposition \ref{CWContribution}]
We have
\begin{align*}
    \mathsf{h}_{X, E, G} &= \sum_{n \geq 0} \frac{1}{n!}\sum_{\sigma \in S_n} \sum_{i \geq 0}(-1)^i \mathrm{Tr}\left(\sigma|_{H^i_c\left(\left(\mathrm{Conf}_n(X) \times \Delta^\circ \right)/G;\QQ\right)}\right) p_1^{k_1(\sigma)}\cdots p_n^{k_n(\sigma)} \\&= \sum_{n \geq 0} \frac{1}{n!}\sum_{\sigma \in S_n} \chi_c\left(\left[\left(\mathrm{Conf}_n(X) \times \Delta^\circ\right)/G\right]^{\sigma}\right) p_1^{k_1(\sigma)}\cdots p_n^{k_n(\sigma)},
\end{align*}
by the Lefschetz fixed-point theorem applied to the one-point compactification of $\left(\mathrm{Conf}_n(X) \times \Delta^\circ\right)/G$, where we set $k_i(\sigma)$ to be the number of cycles of length $i$ in $\sigma$. Now using Lemma \ref{sigmafixedpoints}, we have
\begin{align*}
    \mathsf{h}_{X, E, G} &=-\sum_{n \geq 0} \frac{1}{n!}\sum_{\sigma \in S_n} \frac{(-1)^{|E|}}{|G|} \sum_{g \in G} \mathrm{sgn}(g|_{E}) \cdot \chi_c\left([\mathrm{Conf}_n(X)]^{g^{-1}\sigma}\right) p_1^{k_1(\sigma)}\cdots p_n^{k_n(\sigma)}.
\end{align*}
Now the proof follows that of Gorsky \cite[Theorem 2.5]{gorsky-equivariant}: if we set
\[ X_k(g) := \{x \in X \mid x \text{ has orbit of size }k \text{ under }g \}, \]
and 
\[ \widetilde{X}_k(g) = X_k(g)/(g), \]
then for fixed $\ell_1, \ldots, \ell_n$ such that $\sum_{i = 1}^{n} i \ell_i = n$, we have a map
\[ \coprod_{\substack{{\sigma \in S_n}\\{k_i(\sigma) = \ell_i \forall i}}} [\mathrm{Conf}_n(X)]^{g^{-1}\sigma} \to \prod_{i = 1}^{n} \mathrm{Conf}_{\ell_i}(\widetilde{X}_{i}(g))/S_{\ell_i} \]
which is $n!$-to-$1$, so that
\[\frac{1}{n!} \sum_{\substack{{\sigma \in S_n}\\{k_i(\sigma) = \ell_i \forall i}}} \chi_c\left([\mathrm{Conf}_n(X)]^{g^{-1}\sigma}\right) = \prod_{i = 1}^{n} \frac{\chi_c(\mathrm{Conf}_{\ell_i}(\widetilde{X}_{i}(g)))}{\ell_i!}.  \]
Now the proposition follows from Lemma \ref{ConfigurationSpaceChi}, upon summing over all possible tuples $(\ell_1, \ldots, \ell_n)$.
\end{proof}

Now Theorem \ref{mainthm:FrobChar} is proved by combining Lemma \ref{hgStratificationlemma} with Proposition \ref{CWContribution}.

\begin{remark}\label{remark: generalizedgraphsum}
As explained in Remark \ref{remark: generalizedrepeatedmarkings}, the repeated marking locus in the dual complex $\widetilde{\Theta}_{0, n+N}^{G}(\mu)$ of the inclusion \[\widetilde{\cM}_{0, n + N}^G(\mu) \subset \overline{\cM}_{0, n + N}^G(\mu)\]
is also acyclic, and $\widetilde{\cM}_{0, n + N}^G(\mu)$ is naturally identified with the moduli space of smooth $N$-pointed admissible covers of $\PP^1$ with $\mu$-specified monodromy, together with $n$ distinct marked points on the source curve. 

By the acyclicity of the repeated marking locus, we can write a graph sum formula for the generating function encoding the $S_n$-equivariant weight zero compactly supported Euler characteristics of these moduli spaces. Define
\[ \mathsf{h}^G_{N}(\mu) = \sum_{n \geq 0} \sum_{i = 0}^{2N + 2n - 6} (-1)^i\ch_n(W_0 H^i_c(\widetilde{\cM}_{0, n+N}^G(\mu);\QQ )). \]
By removing the repeated marking locus from the dual complex and emulating the techniques of this section, we obtain the following theorem.
\begin{customthm}{D}\label{mainthm:extendedFrobChar}
We have
\[ \mathsf{h}^G_{N}(\mu) = \sum_{\P \to \C \in \mathrm{Ob}(\Gamma_{0, N}^G(\mu))}  
      \frac{(-1)^{|E_\C|}}{|\Aut(\P \to \C)|} \sum_{\tau \in \Aut(\P \rightarrow \C)} \mathrm{sgn}(\tau|_{E_\C}) \prod_{k \geq 1} (1 + p_k)^{f(\P, \tau, k)}
    \]
    where $E_\C$ is the set of edges of the tree $\C$, 
    $p_k = \sum_{n > 0} x_n^{k} \in \hat{\Lambda}$ is the $k$th power sum symmetric function, and $k \cdot f(\P, \tau, k)$ is given by the compactly supported Euler characteristic of the set of points in $\P$ which have orbit of length $k$, under the action of $\tau$. The first sum is taken over isomorphism classes of objects in $\Gamma_{0, N}^G(\mu)$, which is the category defined in Definition \ref{def:Gamma-category}.
\end{customthm}
 Taking $G = \ZZ/2\ZZ$, $N = 2g + 2$, and $\mu\colon  \{w_1, \ldots, w_N\} \to \ZZ/2\ZZ$ to be the constant function $1$ in Theorem \ref{mainthm:extendedFrobChar}, we obtain the generating function for the $S_n$-equivariant weight zero compactly supported Euler characteristics of the moduli spaces $\widetilde{\cH}_{g, n}$ of $n$-pointed hyperelliptic curves of genus $g$, together with labellings of their Weierstrass points.
\end{remark}

%%%%%%%%%%%%%%%%%%%%%%%%%%%%%%%%%%%%%%%%%%%%%%%%%%%%%%%%%%%%%%%%
\appendix
\section{Calculations for $g\le 7$}\label{sec:calculations}

In this appendix we present the computational data obtained by implementing Theorem \ref{mainthm:FrobChar} on a computer. This was implemented in \emph{Mathematica} using the package IGraph/M \cite{igraphm}. The code for these computations is available at \cite{code}. 

\begin{figure}[h]
    \centering
    \includegraphics{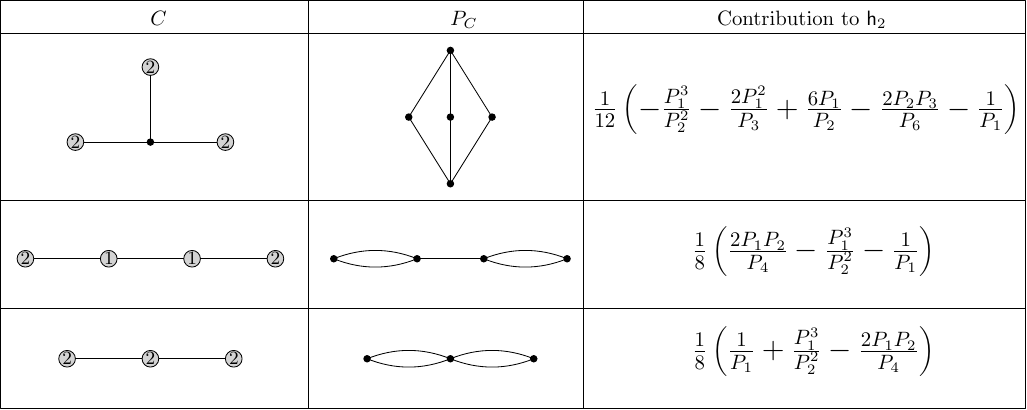}
    \caption{The three trees $C$ in $T_{6}^{< 3}$, their associated covers $P_C$, and the contribution of $P_C \to C$ to $\mathsf{h}_2$ as in Theorem \ref{mainthm:FrobChar}. The generating function $\mathsf{h}_2$ is the sum of the three contributions. Note that the contributions in the second and third rows cancel.}
    \label{fig:h2exmample}
\end{figure}

We compute $\mathsf{h}_g$ explicitly for $2 \le g \le 7$: see Table \ref{table:frobenius_data}. For scale, $\mathsf{h}_5$ is computed as a sum over 96 graphs and takes 8 minutes to compute on a home laptop, while $\mathsf{h}_7$ is computed as a sum over 2789 graphs and takes just under 3 days to compute on a home laptop. Figure \ref{fig:h2exmample} contains the calculation of $\mathsf{h}_2$ as a sum over three graphs; compare with \cite[Example 8.3]{cfgp-sn}.

We extract from this data exponential generating functions for the numerical weight zero compactly supported Euler characteristic by setting $P_1$ to $1+t$ and all other $P_i$ to 1, see Table \ref{table:euler_char_fns}. We display these Euler characteristics for $0 \le n \le 10$ in Table \ref{table:euler_char}.

\begin{table}[]
\renewcommand{\arraystretch}{2}
\begin{tabular}{ | p{0.03\linewidth} | p{0.9\linewidth}| }
\hline
$g$ & $\mathsf{h}_g$                                             \\ \hline
$2$ & 
$
\frac{1}{12} \left(-\frac{P_1^3}{P_2^2}-\frac{2 P_1^2}{P_3}+\frac{6 P_1}{P_2}-\frac{2 P_2 P_3}{P_6}-\frac{1}{P_1}\right)
$                                            \\ \hline
$3$ & $
-\frac{{P_1}^4}{16 {P_2}^3}+\frac{{P_1}^3}{8 {P_2}^2}-\frac{5 {P_1}^2}{16
{P_2}^2}-\frac{{P_1}^2}{8 {P_4}}-\frac{1}{16 {P_1}^2}+\frac{{P_1} {P_2}}{4
{P_4}}+\frac{{P_1}}{2 {P_2}}+\frac{1}{8 {P_1}}-\frac{{P_2}}{8
{P_4}}-\frac{5}{16 {P_2}}
$                                    \\ \hline
$4$ & $
-\frac{9 P_1^5}{160 P_2^4}+\frac{7 P_1^4}{48 P_2^3}-\frac{P_1^3}{8 P_2^2}-\frac{P_1^3}{16 P_2^3}-\frac{P_1^3}{16 P_2 P_4}+\frac{P_1^2}{16 P_2^2}+\frac{P_1^2}{8 P_4}-\frac{P_1^2}{10 P_5}-\frac{P_1}{2 P_2}+\frac{7 P_1}{16 P_2^2}+\frac{P_1}{6 P_3}
-\frac{P_2 P_1}{4 P_4}+\frac{P_1}{8 P_4}+\frac{P_3 P_1}{6 P_6}-\frac{1}{8 P_1}+\frac{7}{48 P_1^2}-\frac{9}{160 P_1^3}+\frac{1}{16 P_2}-\frac{1}{16 P_1 P_2}+\frac{P_2}{8 P_4}-\frac{P_2}{16 P_1 P_4}-\frac{P_2 P_5}{10 P_{10}}
$                                                    \\ \hline
$5$ &
$
 -\frac{11 P_1^6}{192 P_2^5}+\frac{3 P_1^5}{16 P_2^4}-\frac{P_1^4}{4 P_2^3}-\frac{5 P_1^4}{64 P_2^4}-\frac{P_1^4}{16 P_2^2 P_4}+\frac{P_1^3}{8 P_2^2}+\frac{5 P_1^3}{16 P_2^3}+\frac{3 P_1^3}{16 P_2 P_4}-\frac{P_1^2}{4 P_2^2}-\frac{35 P_1^2}{96 P_2^3}-\frac{P_1^2}{12 P_3^2}-\frac{3 P_1^2}{16 P_2 P_4}-\frac{P_2 P_1^2}{8 P_4^2}+\frac{P_1^2}{12 P_6}-\frac{P_4 P_1^2}{8 P_2 P_8}+\frac{3 P_1}{4 P_2^2}-\frac{P_2 P_1}{4 P_4}+\frac{3 P_1}{8 P_4}+\frac{P_2^2 P_1}{4 P_4^2}+\frac{P_4 P_1}{4 P_8}+\frac{1}{8 P_1}-\frac{1}{4 P_1^2}+\frac{3}{16 P_1^3}-\frac{11}{192 P_1^4}-\frac{1}{4 P_2}+\frac{5}{16 P_1 P_2}-\frac{5}{64 P_1^2 P_2}-\frac{35}{96 P_2^2}-\frac{3}{16 P_4}+\frac{3 P_2}{16 P_1 P_4}-\frac{P_2}{16 P_1^2 P_4}-\frac{P_2^2}{8 P_4^2}+\frac{P_2}{12 P_6} -\frac{P_2 P_3^2}{12 P_6^2}-\frac{P_4}{8 P_8}
$                                                    \\ \hline
$6$ & $
 -\frac{227 P_1^7}{3584 P_2^6}+\frac{P_1^6}{4 P_2^5}-\frac{55 P_1^5}{128 P_2^4}-\frac{25 P_1^5}{256 P_2^5}-\frac{9 P_1^5}{128 P_2^3 P_4}+\frac{3 P_1^4}{8 P_2^3}+\frac{7 P_1^4}{16 P_2^4}+\frac{P_1^4}{4 P_2^2 P_4}-\frac{P_1^3}{8 P_2^2}-\frac{13 P_1^3}{16 P_2^3}-\frac{55 P_1^3}{512 P_2^4}-\frac{11 P_1^3}{32 P_2 P_4}-\frac{3 P_1^3}{64 P_2^2 P_4}-\frac{7 P_1^3}{128 P_4^2}-\frac{P_4 P_1^3}{32 P_2^2 P_8}+\frac{7 P_1^2}{8 P_2^2}+\frac{5 P_1^2}{16 P_2^3}+\frac{P_1^2}{4 P_4}+\frac{P_1^2}{8 P_2 P_4}+\frac{P_2 P_1^2}{8 P_4^2}-\frac{P_1^2}{14 P_{7}}-\frac{99 P_1}{64 P_2^2}+\frac{59 P_1}{128 P_2^3}+\frac{P_2 P_1}{4 P_4}-\frac{7 P_1}{8 P_4}+\frac{19 P_1}{64 P_2 P_4}-\frac{9 P_2^2 P_1}{32 P_4^2}+\frac{9 P_2 P_1}{64 P_4^2}-\frac{P_4 P_1}{8 P_8}+\frac{3 P_4 P_1}{16 P_2 P_8}-\frac{1}{8 P_1}+\frac{3}{8 P_1^2}-\frac{55}{128 P_1^3}+\frac{1}{4 P_1^4}-\frac{227}{3584 P_1^5}+\frac{7}{8 P_2}-\frac{13}{16 P_1 P_2}+\frac{7}{16 P_1^2 P_2}-\frac{25}{256 P_1^3 P_2}+\frac{5}{16 P_2^2}-\frac{55}{512 P_1 P_2^2}+\frac{P_2}{4 P_4}+\frac{1}{8 P_4}-\frac{11 P_2}{32 P_1 P_4}-\frac{3}{64 P_1 P_4}+\frac{P_2}{4 P_1^2 P_4}-\frac{9 P_2}{128 P_1^3 P_4}+\frac{P_2^2}{8 P_4^2}-\frac{7 P_2^2}{128 P_1 P_4^2}-\frac{P_4}{32 P_1 P_8}-\frac{P_2 P_{7}}{14 P_{14}}$                                                    \\ \hline
$7$ & $
 -\frac{19 P_1^8}{256 P_2^7}+\frac{351 P_1^7}{1024 P_2^6}-\frac{913 P_1^6}{1280 P_2^5}-\frac{33 P_1^6}{256 P_2^6}-\frac{11 P_1^6}{128 P_2^4 P_4}+\frac{53 P_1^5}{64 P_2^4}+\frac{171 P_1^5}{256 P_2^5}+\frac{185 P_1^5}{512 P_2^3 P_4}-\frac{25 P_1^4}{48 P_2^3}-\frac{389 P_1^4}{256 P_2^4}-\frac{43 P_1^4}{256 P_2^5}-\frac{41 P_1^4}{64 P_2^2 P_4}-\frac{11 P_1^4}{128 P_2^3 P_4}-\frac{P_1^4}{16 P_2 P_4^2}-\frac{P_4 P_1^4}{32 P_2^3 P_8}+\frac{P_1^3}{8 P_2^2}+\frac{15 P_1^3}{8 P_2^3}+\frac{949 P_1^3}{1024 P_2^4}+\frac{9 P_1^3}{16 P_2 P_4}+\frac{29 P_1^3}{64 P_2^2 P_4}+\frac{55 P_1^3}{256 P_4^2}+\frac{5 P_4 P_1^3}{64 P_2^2 P_8}-\frac{23 P_1^2}{16 P_2^2}-\frac{213 P_1^2}{128 P_2^3}-\frac{137 P_1^2}{256 P_2^4}-\frac{3 P_1^2}{8 P_4}-\frac{37 P_1^2}{64 P_2 P_4}-\frac{21 P_1^2}{64 P_2^2 P_4}-\frac{5 P_1^2}{32 P_4^2}-\frac{3 P_2 P_1^2}{64 P_4^2}-\frac{P_2^2 P_1^2}{8 P_4^3}-\frac{P_1^2}{16 P_8}+\frac{P_4 P_1^2}{16 P_2 P_8}-\frac{5 P_4 P_1^2}{32 P_2^2 P_8}+\frac{P_1}{2 P_2}+\frac{57 P_1}{32 P_2^2}+\frac{199 P_1}{128 P_2^3}-\frac{P_1}{6 P_3}+\frac{P_1}{8 P_3^2}+\frac{P_2 P_1}{4 P_4}+\frac{P_1}{4 P_4}+\frac{229 P_1}{256 P_2 P_4}-\frac{5 P_2^2 P_1}{16 P_4^2}+\frac{25 P_2 P_1}{64 P_4^2}+\frac{103 P_2^3 P_1}{384 P_4^3}+\frac{P_1}{10 P_5}-\frac{P_3 P_1}{6 P_6}+\frac{P_3^2 P_1}{8 P_6^2}-\frac{P_4 P_1}{4 P_8}+\frac{5 P_2 P_1}{32 P_8}+\frac{5 P_4 P_1}{16 P_2 P_8}+\frac{P_5 P_1}{10 P_{10}}+\frac{P_6 P_1}{12 P_{12}}+\frac{1}{8 P_1}-\frac{25}{48 P_1^2}+\frac{53}{64 P_1^3}-\frac{913}{1280 P_1^4}+\frac{351}{1024 P_1^5}-\frac{19}{256 P_1^6}-\frac{23}{16 P_2}+\frac{15}{8 P_1 P_2}-\frac{389}{256 P_1^2 P_2}+\frac{171}{256 P_1^3 P_2}-\frac{33}{256 P_1^4 P_2}-\frac{213}{128 P_2^2}+\frac{949}{1024 P_1 P_2^2}-\frac{43}{256 P_1^2 P_2^2}-\frac{137}{256 P_2^3}-\frac{3 P_2}{8 P_4}-\frac{37}{64 P_4}+\frac{9 P_2}{16 P_1 P_4}+\frac{29}{64 P_1 P_4}-\frac{41 P_2}{64 P_1^2 P_4}-\frac{11}{128 P_1^2 P_4}+\frac{185 P_2}{512 P_1^3 P_4}-\frac{11 P_2}{128 P_1^4 P_4}-\frac{21}{64 P_2 P_4} -\frac{5 P_2}{32 P_4^2}-\frac{3 P_2^2}{64 P_4^2}+\frac{55 P_2^2}{256 P_1 P_4^2}-\frac{P_2^2}{16 P_1^2 P_4^2}-\frac{P_2^3}{8 P_4^3}-\frac{P_2}{16 P_8}+\frac{P_4}{16 P_8}+\frac{5 P_4}{64 P_1 P_8}-\frac{P_4}{32 P_1^2 P_8}-\frac{5 P_4}{32 P_2 P_8}
$ \\ \hline
\end{tabular}
\caption{The generating function $\mathsf{h}_g \in \hat{\Lambda}$ for $2 \leq g \le 7$. Here $P_i := 1 + p_i \in \hat{\Lambda}$ is the inhomogeneous power sum.}
\label{table:frobenius_data}
\end{table}

%%%%%%%%

\begin{table}[]
\renewcommand{\arraystretch}{2}
\small
\begin{tabular}{|c|l|}
\hline
$g$ & $\sum_{n \geq 0} \frac{t^n}{n!} \left( \sum_{i = 0}^{4g + 2n - 2} (-1)^i \dim_{\QQ} W_0 H^i_c(\cH_{g, n};\QQ)  \right)$                                          
\\ \hline
$2$ & $- \frac{t^2}{12 (1 + t)}\left(6 + 6 t + t^2\right) $           \\ \hline
$3$ &   $- \frac{t^2}{16 (1 + t)^2}\left(8 + 16 t + 12 t^2 + 4 t^3 + t^4\right) $      \\ \hline
$4$ &   $- \frac{t^2}{480 (1 + t)^3}\left(240 + 720 t + 960 t^2 + 720 t^3 + 386 t^4 + 146 t^5 + 27 t^6\right) $                                                    \\ \hline
$5$ &   $- \frac{t^2}{192 (1 + t)^4}\left(96 + 384 t + 736 t^2 + 864 t^3 + 748 t^4 + 504 t^5 + 246 t^6 + 
 74 t^7 + 11 t^8\right) $  \\ \hline
$6$ &   $- \frac{t^2}{3584 (1 + t)^5}\left(
\begin{aligned}
& 1792 + 8960 t + 22400 t^2 + 35840 t^3 + 43232 t^4 + 41888 t^5 + 
 32096 t^6 + 18272 t^7\\ & + 7268 t^8  + 1828 t^9 + 227 t^{10}
 \end{aligned}
 \right) $  \\ \hline
$7$ & $- \frac{t^2}{15360 (1 + t)^6}\left(
\begin{aligned}
& 7680 + 46080 t + 142080 t^2 + 288000 t^3 + 446720 t^4 + 565760 t^5 + 
 587520 t^6 \\ & + 485120 t^7 + 308936 t^8 + 146832 t^9 + 49551 t^{10} + 
 10695 t^{11} + 1140 t^{12}
 \end{aligned}
 \right) $ 
\\ \hline
\end{tabular}
\caption{The exponential generating functions for numerical weight zero compactly supported Euler characteristics of $\mathcal{H}_{g,n}$.}
\label{table:euler_char_fns}
\end{table}

%%%%%%%%%

\begin{table}[]
\renewcommand{\arraystretch}{2}
\small
\begin{tabular}{|c|l l l l l l l l l l l |}
\hline
 \diagbox{$g$}{$n$}& 0 & 1 & 2 & 3 & 4& 5& 6& 7& 8& 9&  10                 \\ \hline
2& 0& 0& -1& 0& -2& 10& -60& 420& -3360& 30240& -302400 \\
3& 0& 0& -1& 0& -6& 30& -225& 1890& -17640& 181440& -2041200 \\
4& 0& 0& -1& 0& -12& 60& -579& 5586& -59220& 684180& -8557920\\ 
5& 0& 0& -1& 0& -20& 100& -1245& 13230& -157500& 2022300& -27877500\\
6& 0& 0& -1& 0& -30& 150& -2385& 27090& -361080& 5099760& -76856850\\ 
7& 0& 0& -1& 0& -42& 210& -4200& 49980& -745920& 11460960& -187595730\\
 \hline
\end{tabular}
\caption{The weight zero compactly supported Euler characteristic of $\mathcal{H}_{g,n}$ for $2 \leq g \leq 7$, and $0 \leq n \leq 10$.}
\label{table:euler_char}
\end{table}

\clearpage
\bibliographystyle{alpha}
\bibliography{main}

\end{document}

\author[M. Brandt]{Madeline Brandt}
 \address{Department of Mathematics, Vanderbilt University, Nashville, TN 37240}
\email{madeline.v.brandt@vanderbilt.edu}
\author[M. Chan]{Melody Chan}\address{Department of Mathematics, Brown University, Box
1917, Providence, RI 02912}\email{melody\_chan@brown.edu}
\author[S. Kannan]{Siddarth Kannan}\address{Department of Mathematics, Massachusetts Institute of Technology}\email{spkannan@mit.edu}
\title{Erratum to ``On the weight zero compactly supported cohomology of $\mathcal{H}_{g, n}$''}
\maketitle

Let $G$ be a finite abelian group. In \cite{BCK}, we studied the moduli space
\[\Mbar_{0, n}^G(g_1, \ldots, g_n) \]
of \textit{pointed admissible $G$-covers} of $n$-pointed stable curves of genus zero, with monodromies specified by elements $g_i \in G$ satisfying $\sum_i g_i = 0$. This space was introduced by Jarvis--Kimura--Kaufmann \cite{jarvis-kaufmann-kimura-pointed}, and we refer to their work for a detailed definition. In particular, we studied the dual complex of the normal crossings compactification
\begin{equation}\label{eq:normal_crossings}
    \M_{0, n}^G(g_1, \ldots, g_n)  \subset \Mbar_{0, n}^G(g_1, \ldots, g_n)
\end{equation}
where $\M_{0, n}^G(g_1, \ldots, g_n)$ denotes the locus parameterizing covers of smooth genus-zero curves. Proposition 2.4 of loc. cit. purports to calculate the connected components of $\M_{0, n}^G(g_1, \ldots, g_n)$, but it is incorrect as written. The purpose of this note is to prove that Proposition 2.4 is correct in the case $G = \Z/p\Z$ for a prime $p$, which is the case we need for our main application to the moduli space $\mathcal{H}_{g, n}$ of $n$-pointed hyperelliptic curves of genus $g$. We prove this in Proposition \ref{correction} below.

All of the results of \cite{BCK} hold when $G$ is a finite abelian group such that Proposition 2.4 is correct. In particular, the main results Theorem A, Theorem B, and Proposition C are unaffected. In general, our work in \cite{BCK} only determines the dual complex of (\ref{eq:normal_crossings}) for those finite abelian groups $G$ for which Proposition 2.4 holds, so for example Theorem 3.5 and Theorem D only hold for such groups. The groups which have this property can be precisely determined using forthcoming work of Bradford, Pradhan, and Zweifler \cite{BPZ}, which gives a complete description of the connected components of $\M_{0, n}^G(g_1, \ldots, g_n)$ for abelian $G$. See Remark \ref{BPZ} below.
\begin{prop}\label{correction}
    Let $G = \Z/p\Z$ for a prime $p$, and suppose $g_1, \dots, g_n \in G$ satisfy $\sum_{i} g_i = 0$. Then the number of irreducible components of $\M_{0, n}^G(g_1,\ldots,g_n)$ is equal to
    \begin{enumerate}
        \item one, if $g_i \neq 0$ for some $i$;
        \item $p^{n - 1}$ if $g_i = 0$ for all $i$. 
    \end{enumerate}
    As $\M_{0, n}^G(g_1, \ldots, g_n)$ is smooth, the number of irreducible components equals the number of connected components.
\end{prop}
\begin{proof}
    We will use the conventions and notation of \cite[\S 2]{BCK} when working with pointed admissible $G$-covers. Let us first consider the case where $g_i \neq 0$ for some $i$. If in fact $g_i \neq 0$ for \textit{all} $i$, we have an isomorphism \[\M_{0, n}^G(g_1, \ldots, g_n) \cong \M_{0, n} \times BG.\] Therefore to prove the first case of the theorem, we will assume that $g_1 = 0$ and $g_2 \neq 0$. We proceed by induction on $n$. For the base case $n = 3$, we have that $\M_{0, 3}^G(g_1,g_2,g_3)$ is a single point, since the set of lifts of the marked points in the unique admissible $G$-cover $P \to (\mathbb{P}^1, 0, 1, \infty)$ with monodromies specified by $g_i$ is a torsor for $(\prod_{i} G/\langle g_i \rangle)/G$, which is seen to be a single element, since $g_1 + g_2 + g_3 = 0$ implies that $g_2$ and $g_3$ are both generators. For the inductive step, consider the map
    \[ f:\M_{0, n}^G(g_1, \ldots, g_n) \to \M_{0, n - 1}^G(g_2, \ldots, g_n) \]
    which forgets the first marked point and its lift. Note that this map exists because $g_1 = 0$. Observe that $f$ is smooth because it is a morphism between regular Deligne--Mumford stacks of finite type over $\bbC$ whose fibers are regular of constant dimension one: indeed, the fiber of $f$ over a pointed $G$-cover \[[\pi:(P, \tilde{p_2}, \ldots, \tilde{p}_n) \to(C, p_2, \ldots, p_n)] \in \M_{0, n-1}^G(g_2, \ldots, g_n)\] is $P \smallsetminus \pi^{-1}\{p_2, \ldots, p_n \}$. Since $f$ is smooth with an irreducible target and irreducible fibers, we may conclude that $\M_{0, n}^G(g_1, \ldots, g_n)$ is irreducible.\\

    We now turn to the case where $g_i = 0$ for all $i$. We again proceed by induction. When $n = 3$, we see that $\M_{0, 3}^G(0,0,0)$ consists of $p^2$ points. For the inductive step, assume that $\M_{0, n-1}^{G}(0,\ldots,0)$ has $p^{n-2}$ irreducible components. Consider now the auxiliary moduli space 
    \[\M_{0, n-1,*}^{G}(0,\ldots,0) \to \M_{0, n-1}^G(0,\ldots,0).\]
    which parameterizes a point of $\M_{0, n-1}^G(0,\ldots,0)$ together with an additional marked point $p_1$ on the base curve, separate from the other marked points, but \textit{without} a chosen lift of $p_1$ to the source curve. This morphism is also smooth with irreducible fibers, so we can conclude that it induces a bijection on irreducible components. Therefore $\M_{0, n-1,*}^{G}(0,\ldots,0)$ has $p^{n-2}$ irreducible components. Finally, consider the map
    \[ \M_{0, n}^{G}(0,\ldots,0) \to \M_{0, n-1,*}^{G}(0,\ldots,0)\]
    which forgets the lift of $p_1$. This map is a $\Z/p\Z$-torsor, and in fact it is trivial: a section is given by taking the unique lift of $p_1$ which is on the same connected component of the source curve which contains $\tilde{p}_2$. Therefore we have an isomorphism
    \[ \M_{0, n}^{G}(0,\ldots,0) \cong \Z/p\Z \times\M_{0, n-1,*}^{G}(0,\ldots,0), \]
    so $ \M_{0, n}^{G}(0,\ldots,0)$ has $p^{n-1}$ connected components, as desired.
\end{proof}

\begin{remark}\label{BPZ}
    In forthcoming work, Bradford, Pradhan, and Zweifler \cite{BPZ} show that for an arbitrary finite abelian group $G$, there is a natural identification of the connected components of $\M_{0, n}^G(g_1, \ldots, g_n)$ with the orbits of the diagonal $G$-action on the group
    \[ \frac{G^n}{\langle g_i e_i, g_ie_j + g_j e_i \mid i, j \in \{1, \ldots, n\} \rangle}, \]
    where 
    \[ g_i e_j := (0, \ldots, \underbrace{g_i}_{j\text{th entry}}, \ldots, 0) \in G^n. \]
    Their result recovers Proposition \ref{correction} in the case $G = \Z/p\Z$.
\end{remark}

\subsection*{Acknowledgments}
We are grateful to Amy Bradford and  Philip Tosteson for independently notifying us of the error in Proposition 2.4. We also thank Jacob Zweifler for correspondence about the ongoing work \cite{BPZ} and especially for his detailed explanation of the connected components for general abelian groups found in Remark \ref{BPZ}.
\bibliographystyle{alpha}
\bibliography{main_erratum_October2026}